\documentclass{article}

\usepackage{amssymb}

\usepackage{amsthm}
\usepackage{amsmath}

\usepackage{graphics}
\usepackage{float}

\usepackage{hyperref}

\hypersetup{
    colorlinks=true,
    linkcolor=blue
}

\theoremstyle{plain}
\newtheorem{theorem}{Theorem}
\newtheorem{conjecture}{Conjecture}

\title{%
    Conjectures on the Khovanov Homology of Torus Knots, Twist Knots,
    and Legendrian Simple Knots%
}

\author{Vladimir Chernov\hspace{2em}Ryan Maguire}
\date{April 2026}

\begin{document}
    \maketitle
    \begin{abstract}
        \noindent
        A theorem of Kronheimer and Mrowka states that Khovanov homology is
        able to detect the unknot \cite{KronheimerMrowka2011KhovanovUnknot}.
        That is, if a knot has the Khovanov homology of the unknot, then it is
        equivalent to it. Similar results hold for the trefoils
        \cite{BaldwinSivekKhovanovTrefoils} and the figure-eight knot
        \cite{BaldwinDowlinKhovanovFigureEight}. We conjecture that Khovanov
        homology is able to distinguish all torus and twist knots.
        Numerical evidence has been gathered by examining all prime
        knots with 20 or fewer crossings, a total of 2,199,471,680 knots
        (not including mirrors). We found that all knots with the
        same Khovanov polynomial (the Poincar\'{e} polynomial of Khovanov
        homology) as a torus or twist knot are indeed torus or twist knots
        themselves. Since torus knots are known to be Legendrian
        simple, and since all twist knots $K_{m}$ with $m\geq{-3}$ are
        Legendrian simple, this provides evidence for the claim that
        Khovanov homology and Legendrian simplicity may be connected.
        We conjecture that indeed Khovanov homology is able to distinguish
        Legendrian simple knots and use the (conjectured) Legendrian simple
        knots from the Legendrian knot atlas \cite{LegendrianKnotAtlas} to
        test this claim.
        A similar observation was made, and no knots with 20 or fewer crossing
        share their Khovanov polynomial with the knots in the Legendrian
        knot atlas (except for the knots that are a part of this atlas).
        \par\hfill\par
    \end{abstract}
    \section{Legendrian and Transverse Knots and Links}
\begin{figure}
            \centering
            \resizebox{\textwidth}{!}{%
                \includegraphics{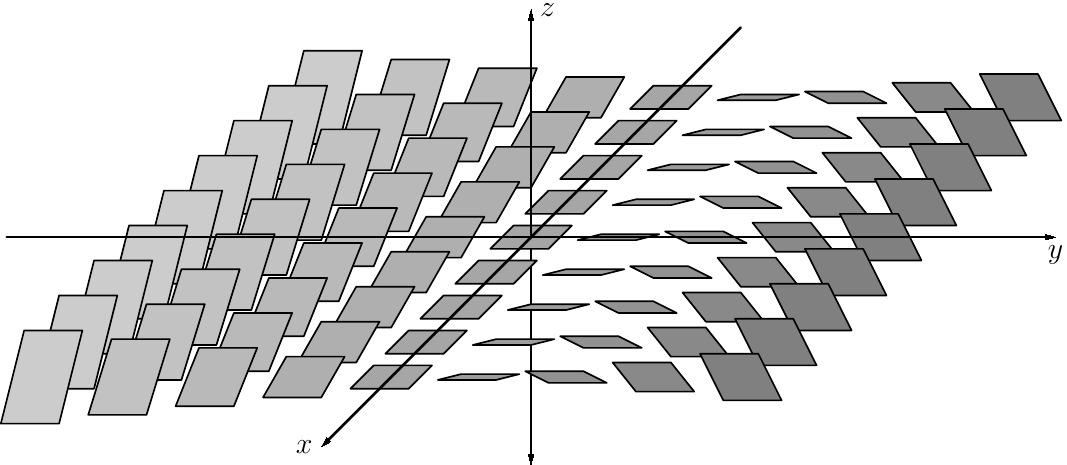}
            }
            \caption{Hyperplane Distribution for $\textrm{d}z-y\,\textrm{d}x$}
            \label{fig:darboux_form_001}
        \end{figure}

        A knot is a smooth embedding of the circle $\mathbb{S}^{1}$ into
        $\mathbb{R}^{3}$. A link is a smooth embedding of $N\in\mathbb{N}$
        disjoint circles in $\mathbb{R}^{3}$. We may impose
        extra structure by considering the standard \textit{contact structure}
        of $\mathbb{R}^{3}$. This is an assignment
        of a plane to every point in $\mathbb{R}^{3}$ such that there is no
        surface $M\subset\mathbb{R}^{3}$ (even infinitesimally) where for every
        $p\in{M}$ the tangent
        plane $T_{p}M$ is given by the plane of the contact structure. In
        $\mathbb{R}^{3}$ this is described by the one-form
        $\textrm{d}z-y\,\textrm{d}x$, the plane at $(x,y,z)$ being spanned by
        the vectors $\partial_{x}+y\partial_{z}$ and $\partial_{y}$. The
        hyperplane distribution is shown in Fig.~\ref{fig:darboux_form_001}.
        While there are no everywhere tangent surfaces, it is possible for a
        curve to be everywhere tangent to this distribution of planes.
        A \textit{Legendrian link} is a link that is everywhere tangent to
        the contact structure. A gentle introduction can be found
        in \cite{JoshuaMSabloffWhatIsLegendrianKnot}.
        \par\hfill\par
        Two Legendrian links are considered to be Legendrian equivalent if
        there is an isotopy $H:L\times[0,1]\rightarrow\mathbb{R}^{3}$ between
        them, $L=\sqcup_{k=0}^{N-1}\mathbb{S}^{1}$, such that for all
        $t\in[0,1]$ the link $H_{t}$ is Legendrian. It is possible for two links
        to be topologically equivalent but not Legendrian equivalent. In the other
        direction, any link can be made Legendrian by an appropriate isotopy
        (see, for example, the introduction of
        \cite{VeraVertessiTransNonSimpleKnots}). The
        two classical invariants of Legendrian links are the Thurston-Bennequin
        and rotation numbers \cite{Ding2006LEGENDRIANHA}.
        A link type is said to be \textit{Legendrian simple} if
        any two Legendrian embeddings of it with the same Thurston-Bennequin
        and rotation numbers are Legendrian equivalent. That is, if the classical
        invariants uniquely classify all Legendrian representations of the knot.
        Certain knot types are known to be Legendrian simple, the first
        discovered being the unknot in 1998
        by Eliashberg and Fraser
        \cite{EliashbergFraserClassificationTopTrivialLegKnots}.
        Torus knots and the figure eight knot are also known to be Legendrian
        simple
        \cite{EtnyreHondaContactTopologyI}, and in
        \cite{EtnyreEtAlLegendrianAndTransverseTwistKnots} Etnyre, Ng, and
        Vertesi classify which twist knots are Legendrian simple. In particular,
        the $K_{m}$ twist knot is Legendrian simple if and only if
        $m\geq{-3}$.
        \par\hfill\par
        Not every knot is Legendrian simple, the mirror of the $K_{3}$ twist
        knot (also known as the $5_{2}$ knot, which is isotopic to the
        $K_{-4}$ twist knot) being an example discovered by Chekanov
        \cite{ChekanovDifAlgOfLegLinks}, and independently by
        Eliashberg \cite{EliashbergInvariantsInContactTopology}.
        More examples can be found in the
        Legendrian knot atlas \cite{LegendrianKnotAtlas}.
        \par\hfill\par
        Transverse links are links that are everywhere transverse to the
        contact structure. That is, at every point $p$ on the link the velocity
        vector and the contact hyperplane at $p$ span $\mathbb{R}^{3}$.
        Any Legendrian
        link can be made transverse by a small perturbation in the direction
        normal to the given plane in the contact structure.
        Two transverse links are transversely equivalent if there is an
        isotopy $H:L\times[0,1]\rightarrow\mathbb{R}^{3}$ such that $H_{t}$ is a
        transverse link for all $t\in[0,1]$. The Bennequin number of a
        transverse knot is defined by the \textit{algebraic crossing number}
        $w(K)$ (more commonly known as the \textit{writhe}) and the\
        \textit{braid index} $n(K)$. It is
        \begin{equation}
            \beta(K)=w(K)-n(K).
        \end{equation}
        It is not an invariant of topological knots, but is an invariant under
        transverse equivalence. This is essentially the self-linking number of
        the corresponding topological knot with respect to the natural framing
        coming from the trivialization of the contact structure.
        Similar to Legendrian simple, we define a knot
        (or link) type to be transversely simple if all of its transverse
        representations are uniquely determined by their Bennequin number
        (see \cite{BirmanWrinkleTransversallySimpleKnots}) and by whether its
        velocity vectors point into the half space where the contact structure
        is positive or not. A paper by Etynre, Ng, and Vertesi
        \cite{EtnyreEtAlLegendrianAndTransverseTwistKnots}
        classifies when twist knots are transversely simple. In
        particular, infinitely many such knots are transversely simple%
        \footnote{%
            The twist knot $K_{m}$ is transversely simple for $m\geq{-2}$
            or $m$ odd \cite{EtnyreEtAlLegendrianAndTransverseTwistKnots}.
        }
        giving
        us a family of knots to test conjectures with. It is also true that
        infinite families of non-Legendrian simple and non-transversely simple
        knots exist. See, for example, the works of Etnyre and Honda
        \cite{EtnyreHondaCabling} and Birman and Menasco
        \cite{BirmanMenasco2006}. More examples can be found in
        \cite{Foldvari2019legnonsimple}.
        \par\hfill\par
        There are several common ways of representing topological knots, the
        three used in our computations are extended Gauss code, planar diagram
        code (PD code), and Dowker-Thistlewaite code (DT code). Given a knot
        diagram with $N$ crossings, the (unsigned) Gauss code is a string with $2N$
        characters (extended Gauss code is $3N$),
        PD code is a string that is $4N$ long, and DT code is $N$
        characters long. Each has its benefits. Extended Gauss code
        can distinguish mirrors, DT code cannot
        (see \cite{DowkerThistlethwaiteDTCode}).
        PD code is perhaps the easiest to
        reconstruct the knot diagram, and DT code is the shortest. Because of
        this we will present our examples in DT code. To obtain the DT code of
        a knot diagram, place your finger on the knot and \textit{walk} along
        the diagram, labelling the crossings. When you get back to your starting
        point each crossing will have two numbers associated with it. It is
        not difficult to see that each crossing will have exactly one odd number
        and one even number.%
        \footnote{%
            For a \textit{classical} knot,
            one that embeds into $\mathbb{R}^{3}$.
            This claim is not true for virtual knots.
        }
        For each even number, if that number was associated
        with an \textit{over} crossing (that is, your finger ran over the
        crossing as you were labelling it), place a minus sign in front. Write
        out the pairs of integers as $(1,a_{1})$, $(3,a_{2})$, $\dots$,
        $(2n-1,a_{n})$. The DT code is the list $a_{1},a_{2},\dots,a_{n}$.
        See \cite{KatlasDTCode} for several examples.
        \par\hfill\par
        It is possible to go from DT code to unsigned
        Gauss code (i.e. the usual Gauss code, and not the extended Gauss code)
        and back \cite{KatlasDTCode}. For certain computations,
        like the Alexander polynomial which
        is mirror insensitive, DT code is easiest since it is the shortest.
        For things like the Jones polynomial and Khovanov homology, invariants
        that distinguish mirrors, extended Gauss code is a must.
        \par\hfill\par
        A note on notation. The tables on the following pages identify knots
        using their DT codes. The \texttt{regina} library \cite{regina} uses a
        clever trick to efficiently express a knots DT code
        for all knots of up to 26 crossings, and we'll
        adopt this notation for the remainder of the paper. The DT code of an
        $n$ crossing knot is a string of length $n$ of signed even integers up
        to $2n$. Since all entries are even, we may divide by two and obtain
        a string of length $n$ of signed integers up to $n$. So long as $n$ is
        not larger than 26, we may replace integers with letters. We then use
        lower case letters for positive integers, and upper case for negative.
        For example, the trefoil $T(2,\,3)$ becomes:
        \begin{equation}
            4,\,6,\,2
            \Rightarrow
            2,\,3,\,1
            \Rightarrow
            \texttt{bca}
        \end{equation}
        The knot $10_{132}$, which is the 10 crossing knot in
        Tab.~\ref{table:matching_twist_knots} with the same Jones polynomial
        as $T(2,\,5)$, becomes:
        \begin{align}
            8,\,6,\,18,\,2,\,-12,\,-16,\,20,\,-10,\,4,\,14
            &\Rightarrow
            4,\,3,\,9,\,1,\,-6,\,-8,\,10,\,-5,\,2,\,7\\
            &\Rightarrow
            \texttt{dciaFHjEbg}
        \end{align}
    \section{Khovanov Homology}
        The Khovanov homology of a link is a powerful, if computationally
        expensive%
        \footnote{%
            The na\"{i}ve algorithm is exponential in the number of
            crossings. Improvements by Bar-Natan \cite{BarNatan2006FASTKH}
            have sped up computations but no polynomial-time algorithm is
            known at the time of this writing.
        },
        invariant first introduced by Mikhail Khovanov
        \cite{Khovanov1999CatJonesPoly} (see also \cite{BarNatanKhovanovJones}
        for an excellent introduction). It is closely related to the Jones
        polynomial, but able to distinguish many more knots and links. The
        homology groups $KH^{r}(L)$ of a link (or knot) $L$ are the direct sum
        of homogeneous components $KH_{\ell}^{r}(L)$ and the
        \textit{Khovanov Polynomial} (see \cite{KatlasKhoHo}) is given by
        \begin{equation}
            Kh(L)(q,\,t)
            =
            \sum_{r,\ell}t^{r}q^{\ell}\textrm{dim}\big(KH_{\ell}^{r}(L)\big).
        \end{equation}
        The Jones polynomial of $L$ is recovered via
        \begin{equation}
            J(L)(q)
            =
            Kh(L)(q,\,-1).
        \end{equation}
        Khovanov homology is not a perfect invariant
        \cite{Watson2007KnotsWI}. That is, there are
        distinct knots with the same Khovanov homology, but it is a powerful
        invariant and is capable of detecting certain knot types.
        \begin{theorem}[Kronheimer and Mrowka, 2011]
            If a knot $K$ has the same Khovanov homology as the unknot, then $K$
            is equivalent to the unknot.
        \end{theorem}
        The unknotting problem asks one to determine if a given knot diagram is
        equivalent to the unknot. Khovanov homology is a powerful enough tool
        to accomplish this task. The Khovanov polynomial is a generalization of
        the Jones polynomial and it has been conjectured that if a
        knot has the same Jones polynomial as the unknot, then that knot is
        equivalent to the unknot. At the time of this writing it has not been
        proven, but there is evidence for and against the claim.
        Thistlewaite found links with the same Jones polynomial as the unlink
        \cite{Thistlethwaite2001LINKSWT}, and there is a 3-crossing virtual
        knot that has the same Jones polynomial as the unknot. For the claim,
        all knots of up to 24 crossings are either the unknot, or have a
        Jones polynomial different from the unknot
        \cite{VerificationUnknotJonesConjUpTo24}.
        \par\hfill\par
        In recent years it has been discovered that Khovanov homology can detect
        a few other knot types. In particular, both the left and right handed
        trefoils, as well as the figure eight.
        \begin{theorem}[Baldwin and Sivek, 2022]
            If a knot $K$ has the same Khovanov homology as either of the
            trefoils, then $K$ is equivalent to one of them.
        \end{theorem}
        \begin{theorem}[Baldwin, Dowlin, Levine, Lidman, and Sazdanovic, 2021]
            If a knot $K$ has the same Khovanov homology as the figure-eight
            knot, then $K$ is equivalent to it.
        \end{theorem}
        See \cite{BaldwinSivekKhovanovTrefoils} and
        \cite{BaldwinDowlinKhovanovFigureEight}, respectively.
    \section{Conjectures on Khovanov Homology}
        Khovanov homology is capable of detecting the
        unknot and both the right and left trefoils. The Khovanov homology with
        coefficients in $\mathbb{Z}/2\mathbb{Z}$ is also capable of detecting
        the cinquefoil knot \cite{BaldwinYingSivekCinquefoilKhovanov},
        which is the $T(5,2)$ torus knot. The Jones
        polynomial, on the other hand, is not capable of detecting the
        $T(5,2)$ torus knot since the $10_{132}$ knot yields the same
        polynomial (see Tab.~\ref{table:matching_torus_knots}).
        We make the following conjecture.
        \begin{conjecture}
            If $K$ is a torus knot, $K=T(p,\,q)$ for co-prime $p$ and $q$,
            then the Khovanov homology of $K$ distinguishes $K$. That is, if
            $K^{\prime}$ is another knot and
            $Kh(K)=Kh(K^{\prime})$, then $K$ and $K^{\prime}$ are equivalent.
        \end{conjecture}
        Numerical evidence has been tallied by comparing the Khovanov homologies
        of torus knots  against all prime knots up to 20 crossings, a total of
        2,199,471,680 knots
        \cite{Burton2020TheN3, Thistlethwaite20CrossingKnots}.
        There are four
        non-torus knots that have the same Jones polynomial as a torus knot
        ($T(2,5)$ matches 10 and 17 crossing knots, $T(2,7)$
        matches a 12 crossing knot, and $T(2,11)$ matches a 14 crossing knot)
        so we cannot generalize the Jones unknot conjecture. Nevertheless, in
        all cases the Khovanov homologies were different
        (see Numerical Results).
        \par\hfill\par
        We also performed the computation with twist knots, and found that
        there are no knots with 20 or fewer crossings that have the same
        Khovanov polynomial as a twist knot (other than twist knots themselves).
        In light of this, we make the following conjecture.
        \begin{conjecture}
            If $K$ is a twist knot, $K=K_{m}$ for some $m\in\mathbb{Z}$,
            then the Khovanov homology of $K$ distinguishes $K$.
        \end{conjecture}
        Since the torus knots are known to be Legendrian simple
        \cite{EtnyreHondaContactTopologyI}, and since the twist knot $K_{m}$
        is Legendrian simple for all $m\geq{-3}$
        \cite{EtnyreEtAlLegendrianAndTransverseTwistKnots},
        the fact that both of these
        families appear to be detected by their Khovanov homologies suggests
        that Khovanov homology may detect Legendrian simplicity. We repeated
        our calculations using the 26 conjecturally Legendrian simple knots
        from the Legendrian knot atlas \cite{LegendrianKnotAtlas},%
        \footnote{%
            There are several ways to count distinct knots since a knot,
            its mirror, its inverse, and its inverse mirror may all be distinct.
            We do not distinguish a knot from its inverse in our counting,
            but do treat a knot and its mirror with more care.
        }
        and again found that there were no knots with 20 or fewer crossings
        that shared their Khovanov homology with the knots in this table.
        Hence, we conjecture the following.
        \begin{conjecture}
            If $K$ is a Legendrian simple knot,
            then the Khovanov homology of $K$ distinguishes $K$.
        \end{conjecture}
        We explore each of these three cases later in the Numerical Results
        sections.
        \par\hfill\par
        The computations were done as follows. There are libraries in Python,
        Sage, and C++ for working with knot polynomials. In particular, we used
        Regina \cite{regina}, SnapPy \cite{SnapPy}, the Sage knot library
        \cite{sage}, and our own ever-growing C library
        \cite{MaguireLibtmpl}.
        \par\hfill\par
        We first gather the torus knots that could potentially have the same
        Jones polynomial as some other knot based on the bounds of the degree.
        The Jones polynomials of torus knots were computed using the formula
        \cite{JonesJonesPolyForDummies}
        \begin{equation}
            \label{eqn:jones_poly_torus}%
            J(T(m,n))(q)=
            q
            ^{(m-1)(n-1)/2}
            \frac{1-q^{m+1}-q^{n+1}+q^{m+n}}{1-q^{2}}.
        \end{equation}
        The Kauffman bracket, of which the Jones polynomial is a simple
        normalization \cite{BarNatanKhovanovJones},
        of a knot $K$ with $N$ crossings has the formula
        \begin{equation}
            \langle{K}\rangle
            =
            \sum_{n=0}^{2^{N}-1}(-q)^{w(n)}(q+q^{-1})^{c(n)},
        \end{equation}
        where $w(n)$ is the Hamming weight of $n$, the number of 1's in the
        binary representation of $n$, and $c(n)$ is the number of disjoint
        cycles that result from the complete smoothing of $K$ corresponding to
        $0\leq{n}\leq{2}^{N}-1$. With this the degree of the bracket polynomial
        is bounded by $2N+1$ ($w(n)$ is bounded by $N$, and $c(n)$ is bounded
        by $N+1$). The degree of the normalization is bounded by $2N$
        \cite{BarNatanKhovanovJones}, so the
        Jones polynomial of a knot with $N$ crossings has degree at most
        $4N+1$.
        \par\hfill\par
        Since we are looking at knots with up to $20$ crossings, we collect the
        co-prime pairs $(m,\,n)$ with $1<m<n$ such that the degree of
        $J\big(T(m,\,n))(q)$ is less than $4\cdot{20}+1=81$.
        Since the Jones polynomial of a mirror
        can be computed by substituting $q\mapsto{q}^{-1}$
        \cite{JonesPolyJones}, we need not look at
        negative values. The Khovanov polynomial makes a similar change,
        $(q,\,t)\mapsto(q^{-1},\,t^{-1})$
        \cite{WatsonKhovanovHomology2017}.
        \par\hfill\par
        Using any of the aforementioned libraries, the Jones polynomials of all
        prime knots up to 20 crossings were computed and compared against this
        table of torus knot Jones polynomials (Eqn.~\ref{eqn:jones_poly_torus}).
        If a match was found the regina library was used to determine if the
        knots were actually identical. That is, if the knot whose Jones
        polynomial was being compared against the torus knots was indeed a
        torus knot itself. If the knots were distinct, this knot was saved in a
        text file for later examination. At the end of the computation 4
        non-torus knots had the same Jones polynomial as a torus knot
        (see Tab.~\ref{table:matching_torus_knots}).
        Since the Khovanov polynomial contains the Jones polynomial in it
        (recall $J(L)(q)=Kh(L)(q,\,-1)$) the only possible non-torus knots with
        the same Khovanov homology as a torus knot were these 4.
        \par\hfill\par
        Using the Java library JavaKh\footnote{%
            Thanks must be paid to Nikolay Pultsin who made edits to
            JavaKh-v2 so that it may run on a GNU/Linux machine using
            OpenJDK 17.
        }
        we found that these four knots with the same Jones polynomials as some
        torus knot all had different Khovanov homologies. Thus, we have the
        following claim.
        \begin{theorem}
            If a prime knot $K$ has less than or equal to 20 crossings and has
            the Khovanov homology of a torus knot $T$,
            then $K$ is equivalent to $T$.
        \end{theorem}
        A similar search through the twist knots yielded more results.
        The Jones polynomials of the twist knots are known, with the formula:
        \begin{equation}
            J(K_{n})(q)=
            \begin{cases}
                (1+q^{-2}+q^{-n}+q^{-n-3})/(1+q),&n\textrm{ odd.}\\
                (1+q-q^{3-n}+q^{-n})/(1+q),&n\textrm{ even.}
            \end{cases}
        \end{equation}
        A search through all prime knots with up to 20 crossings against twist
        knots provided 11 matches for the Jones polynomial
        (see Tab.~\ref{table:matching_twist_knots}), but
        none for Khovanov homology.
        \par\hfill\par
        Since infinitely many of the twist knots are transversely simple, it is
        tempting to conjecture that transverse simplicity may be related to
        Khovanov homology.
        Unfortunately, in \cite{Petkova02012026} Petkova and Schwartz provide
        numerical evidence that transverse non-simplicity is not very common.
        Through their calculations, they conjecture that, of the 604 knots of
        arc index 10 that were examined, 135 of them are transversely
        non-simple, or about 22\%, meaning roughly 78\% of these knots are
        transversely simple. Our calculations with Khovanov homology show that
        for knots with 17 crossings or less, roughly 32\% of knots have
        Khovanov homologies that uniquely distinguish them among other knots
        with 17 crossings (or fewer) (see Tab.~\ref{tab:kho_stat}).
        Thus the conjecture in \cite{Petkova02012026} and
        the conjecture that Khovanov homology detects transverse simplicity
        are incompatible; either many of the conjectured
        transversely simple knots in \cite{Petkova02012026} are actually
        transversely non-simple, or Khovanov homology does not detect
        transversely simplicity. Since only twist knots were used in our
        search, it is likely that the latter is the false conjecture.
        Nevertheless, experimenting with twist knots does show that, for up to
        20 crossings, twist knots do not share their Khovanov homology with any
        other knots. We examine this further in the Numerical Results section.
        \par\hfill\par
        We also performed our computations on the conjecturally Legendrian
        simple knots in \cite{LegendrianKnotAtlas}. Like the twist knots,
        no matches were found for the Khovanov polynomial.
        \par\hfill\par
        A similar search using Knot Floer Homology (KFH), a homology theory
        first introduced by Peter Ozsv\'{a}th and Zolt\'{a}n Szab\'{o}
        \cite{OzsvathSzaboHolomorphicDisks2004}, using the Alexander polynomial
        was performed. This search found that there are several distinct knots
        with the same Knot Floer Homology as a Legendrian simple knot,
        including torus and twist knots. Steven Sivek pointed
        out that the pretzel knots $P(-3,\,3,\,2n+1)$ all have the same KFH,
        meaning the $6_{1}$ twist knot matches the KFH of the $9_{46}$
        knot in the Rolfsen table.
        Matthew Hedden also informed us that the $T(4,\,3)$ knot and the
        $(2,3)$ cable of the trefoil also have matching KFH.
    \section{Numerical Results}
        \subsection{Torus Knots}
        \begin{table}[H]
            \centering
            \resizebox{\textwidth}{!}{%
                \begin{tabular}{| l | l | l |}
                    \hline
                        Torus Knot&
                        Non-Torus Knot&
                        Jones Polynomial\\
                    \hline
                        $T(2,\,5)$&
                        \texttt{dciaFHjEbg}&
                        $-q^{14}+q^{12}-q^{10}+q^{8}+q^{4}$\\
                    \hline
                        $T(2,\,7)$&
                        \texttt{fJGkHlICEABd}&
                        $-q^{20}+q^{18}-q^{16}+q^{14}-q^{12}+q^{10}+q^{6}$\\
                    \hline
                        $T(2,\,11)$&
                        \texttt{gHlImJnKBDFAce}&
                        $-q^{32}+q^{30}-q^{28}+q^{26}-q^{24}+q^{22}-q^{20}+q^{18}-q^{16}+q^{14}+q^{10}$\\
                    \hline
                        $T(2,\,5)$&
                        \texttt{iNHlPJqCoKFmdABgE}&
                        $-q^{14}+q^{12}-q^{10}+q^{8}+q^{4}$\\
                    \hline
                \end{tabular}%
            }
            \caption{Knots whose Jones polynomial matches that of a Torus Knot}
            \label{table:matching_torus_knots}
        \end{table}
        In Tab.~\ref{table:matching_torus_knots}
        we have found four non-torus knots with the same
        Jones polynomial as some torus knot. In all four cases the Khovanov
        polynomials differ. The tables for these polynomials are given
        in Appendix A.
        The coefficient of $q^{\ell}t^{r}$ is given by the corresponding entry.
        Empty equates to zero.
        \par\hfill\par
        The cinquefoil $T(2,\,5)$ has the same Jones polynomial as two
        non-torus knots. As shown in
        Tabs.~\ref{table:t_2_5_kho}-\ref{table:iNHlPJqCoKFmdABgE_kho}, the
        Khovanov polynomials differ
        (this had already been known, see \cite{KatlasKhoHo}).
        Note the columns for the Euler characteristic $\chi$ are identical,
        indicating matching Jones polynomials.
        \par\hfill\par
        The $T(7,\,2)$ knot, also the $7_{1}$ knot, and occasionally called the
        septafoil, has the same Jones polynomial as \texttt{fJGkHlICEABd}. The
        Khovanov polynomials are distinct
        (Tabs.~\ref{table:t_7_2_kho}-\ref{table:fJGkHlICEABd_kho}).
        Lastly, the $T(11,\,2)$ torus knot has the same Jones polynomial as the
        14 crossing knot \texttt{gHlImJnKBDFAce}. Once again the Khovanov
        polynomials differ.
        \subsection{Twist Knots}
        \begin{table}[H]
            \centering
            \resizebox{\textwidth}{!}{%
                \begin{tabular}{| l | l | l |}
                    \hline
                        Twist Knot&
                        Non-Twist Knot&
                        Jones Polynomial\\
                    \hline
                        $K_{2}$&
                        \texttt{eikGbHJCaFd}&
                        $q^{4}-q^{2}+1-q^{-2}+q^{-4}$\\
                    \hline
                        $K_{3}$&
                        \texttt{dgikFHaEjbc}&
                        $-q^{12}+q^{10}-q^{8}+2q^{6}-q^{4}+q^{2}$\\
                    \hline
                        $K_{3}$&
                        \texttt{gfJKHlaIEBCD}&
                        $-q^{12}+q^{10}-q^{8}+2q^{6}-q^{4}+q^{2}$\\
                    \hline
                        $K_{3}$&
                        \texttt{hGJaMlCdEKBfI}&
                        $-q^{12}+q^{10}-q^{8}+2q^{6}-q^{4}+q^{2}$\\
                    \hline
                        $K_{5}$&
                        \texttt{bhDGijCkaef}&
                        $-q^{16}+q^{14}-q^{12}+2q^{10}-2q^{8}+2q^{6}-q^{4}+q^{2}$\\
                    \hline
                        $K_{6}$&
                        \texttt{cefIgbajkDh}&
                        $q^{12}-q^{10}+q^{8}-2q^{6}+2q^{4}-2q^{2}+2-q^{-2}+q^{-4}$\\
                    \hline
                        $K_{6}$&
                        \texttt{femIbaJKLCGHd}&
                        $q^{12}-q^{10}+q^{8}-2q^{6}+2q^{4}-2q^{2}+2-q^{-2}+q^{-4}$\\
                    \hline
                        $K_{6}$&
                        \texttt{jpIFNMrClqOhkEDabg}&
                        $q^{12}-q^{10}+q^{8}-2q^{6}+2q^{4}-2q^{2}+2-q^{-2}+q^{-4}$\\
                    \hline
                        $K_{7}$&
                        \texttt{cgjFHIaDEkb}&
                        $-q^{20}+q^{18}-q^{16}+2q^{14}-2q^{12}+2q^{10}-2q^{8}+2q^{6}-q^{4}+q^{2}$\\
                    \hline
                        $K_{8}$&
                        \texttt{knIHoBjCDQrMPaeLgF}&
                        $q^{16}-q^{14}+q^{12}-2q^{10}+2q^{8}-2q^{6}+2q^{4}-2q^{2}+2-q^{-2}+q^{-4}$\\
                    \hline
                        $K_{9}$&
                        jopIFMrDlqNhkEabcg&
                        $-q^{24}+q^{22}-q^{20}+2q^{18}-2q^{16}+2q^{14}-2q^{12}+2q^{10}-2q^{8}+2q^{6}-q^{4}+q^{2}$\\
                    \hline
                \end{tabular}%
            }
            \caption{Knots whose Jones polynomial matches that of a Twist Knot}
            \label{table:matching_twist_knots}
        \end{table}
        The twist knots have a few more matching non-twists knots for the
        Jones polynomial (Tab.~\ref{table:matching_twist_knots}).
        All tables are provided in Appendix B.
        The first match is the
        figure-eight knot, $K_{2}$, which has the same Jones polynomial as
        K11n19 from the Hoste-Thistlewaite table.
        Again, this had already been known
        (see \cite{KatlasFigureEight} and \cite{KatlasK11n19}).
        The Khovanov polynomials are
        given in Tabs.~\ref{table:m_2_kho}-\ref{table:eikGbHJCaFd_kho}.
        \par\hfill\par
        The $K_{-4}$ twist knot, which is $5_{2}$ on the Rolfsen table,
        has the same Jones polynomial as (at least) three other knots. In each
        case the Khovanov polynomial distinguishes it
        (Tabs.~\ref{table:m_3_kho}-\ref{table:hGJaMlCdEKBfI_kho}).
        The $K_{5}$ twist knot, which is $7_{2}$ on the Rolfsen table, has the
        same Jones polynomial as \texttt{bhDGijCkaef}. The Khovanov polynomials
        are given in Tabs.~\ref{table:m_5_kho} and \ref{table:bhDGijCkaef_kho}.
        \par\hfill\par
        Our original experiment concluded at 17 crossings. We expanded
        this to 19 in 2023, and 20 crossings in 2025%
        \footnote{
            We express our gratitude towards Morwen Thistlethwaite for
            sharing his tabulations with us, see
            \cite{Thistlethwaite20CrossingKnots}.
        },
        and in doing so we have added a few new results,
        the first of which comes with
        $K_{6}$. An 18 crossing knot was found
        that has the same Jones polynomial, but the Khovanov polynomials
        differ. This is shown in
        Tabs.~\ref{table:m_6_kho} and \ref{table:jpIFNMrClqOhkEDabg_kho}.
        $K_{6}$ also shares its Jones polynomial with an 11 and a 13 crossing
        knot. The Khovanov polynomials are given in
        Tabs.~\ref{table:cefIgbajkDh_kho} and
        \ref{table:femIbaJKLCGHd_kho}. Something to note is that the
        the Khovanov polynomials of \texttt{cefIgbajkDh} and
        \texttt{femIbaJKLCGHd} are identical, indicating that the
        Khovanov polynomial is not a perfect invariant.
        \par\hfill\par
        Lastly, the twist knots $K_{7}$, $K_{8}$, and $K_{9}$ have the same
        Jones polynomial of at least one other knot. As before, the Khovanov
        polynomials differ in each case. We may now make the following claim.
        \begin{theorem}
            If a prime knot $K$ has less than or equal to 20 crossings and has
            the Khovanov homology of a twist knot $T$,
            then $K$ is equivalent to $T$.
        \end{theorem}
        An interesting thing to note is that adding 19 and 20 crossing knots
        into the search did not change much. There were no twist knots
        (or torus knots) that have the same Jones polynomial as a 19 or 20
        crossing knot. It seems to be a rare phenomenon for torus and twist
        knots to share their Jones polynomial with low crossing knots.
        \subsection{Conjectured Legendrian Simple Knots}
        Ng presents several knots in \cite{LegendrianKnotAtlas} that are
        \textit{possibly} Legendrian simple but not confirmed. For completeness
        we used our scripts on these knots as well and found several knots with
        identical Jones polynomials, but again all had differing Khovanov
        polynomials. Extending our search from 19 to 20 crossings only
        introduced one Jones collision:
        the knot \texttt{bfINqJaKLTGHPRsCeMDo} has the same Jones polynomial as
        the \texttt{K11n118} knot, but the Khovanov homologies differ.
        All such matches are listed int
        Tab.~\ref{fig:conj_leg_simp_knots} ($m$ indicates the knots mirror).
        Note this table also includes knots that are conjectured to be
        transversely simple.
        The tables for the Khovanov polynomials are given in Appendix C.
        \begin{table}
            \centering
            \resizebox{\textwidth}{!}{%
                \begin{tabular}{| l | l | l |}
                    \hline
                    Ng Knot&Matching Knot&Jones Polynomial\\
                    \hline
                    $m(6_{2})$&\texttt{glfoJcbKMNDaHIe}&$q^{10}-2q^{8}+2q^{6}-2q^{4}+2q^{2}-1+q^{-2}$\\
                    \hline
                    $m(6_{2})$&\texttt{hknEGmDbJLaIfc}&$q^{10}-2q^{8}+2q^{6}-2q^{4}+2q^{2}-1+q^{-2}$\\
                    \hline
                    $m(6_{2})$&\texttt{gKHlmIdJCEABf}&$q^{10}-2q^{8}+2q^{6}-2q^{4}+2q^{2}-1+q^{-2}$\\
                    \hline
                    $m(6_{2})$&\texttt{ehkjmGIaFlcbd}&$q^{10}-2q^{8}+2q^{6}-2q^{4}+2q^{2}-1+q^{-2}$\\
                    \hline
                    $m(7_{3})$&\texttt{hgelkIbaJFcd}&$-q^{18}+q^{16}-2q^{14}+3q^{12}-2q^{10}+2q^{8}-q^{6}+q^{4}$\\
                    \hline
                    $m(7_{4})$&\texttt{gfkHlbjIDaec}&$-q^{16}+q^{14}-2q^{12}+3q^{10}-2q^{8}+3q^{6}-2q^{4}+q^{2}$\\
                    \hline
                    $m(9_{48})$&\texttt{gnoqKDjIMrpEaHblfc}&$q^{2}-3+4q^{-2}-4q^{-4}+6q^{-6}-4q^{-8}+3q^{-10}-2q^{-12}$\\
                    \hline
                    $m(9_{49})$&\texttt{lFKJIOAEnDCpBhmG}&$q^{-4}-2q^{-6}+4q^{-8}-4q^{-10}+5q^{-12}-4q^{-14}+3q^{-16}-2q^{-18}$\\
                    \hline
                    $m(10_{128})$&\texttt{eHPNqGJlBFoiaDCkM}&$-q^{20}+q^{18}-2q^{16}+2q^{14}-q^{12}+2q^{10}-q^{8}+q^{6}$\\
                    \hline
                    $m(10_{128})$&\texttt{edjkaGIlFbch}&$-q^{20}+q^{18}-2q^{16}+2q^{14}-q^{12}+2q^{10}-q^{8}+q^{6}$\\
                    \hline
                    $m(10_{136})$&\texttt{igDKHJaEbFC}&$q^{6}-2q^{4}+2q^{2}-2+3q^{-2}-2q^{-4}+2q^{-6}-q^{-8}$\\
                    \hline
                    $10_{145}$&\texttt{eoHKqGJnCFmPDibaL}&$-q^{20}+q^{18}-q^{16}+q^{14}+q^{4}$\\
                    \hline
                    $10_{145}$&\texttt{kNJIpHLFECoMGABd}&$-q^{20}+q^{18}-q^{16}+q^{14}+q^{4}$\\
                    \hline
                    $10_{161}$&\texttt{hOqrljsnMeipFAgkbcd}&$-q^{22}+q^{20}-q^{18}+q^{16}-q^{14}+q^{12}+q^{6}$\\
                    \hline
                    \texttt{K11n118}
                    &
                    \texttt{bfINqJaKLTGHPRsCeMDo}
                    &
                    $-q^{18}+2q^{16}-3q^{14}+4q^{12}-4q^{10}+3q^{8}-2q^{6}+2q^{4}$\\
                    \hline
                \end{tabular}%
            }
            \caption{Conjectured Legendrian Simple Knots}
            \label{fig:conj_leg_simp_knots}
        \end{table}
    \section{Jones and Khovanov Polynomial Databases}
        Using the \texttt{regina} and \texttt{JavaKh} libraries, the Jones
        polynomials of all knots with less than or equal to 20 crossings were
        tabulated, and the Khovanov polynomials of all knots up to 17 crossings
        were computed.%
        \footnote{%
            This data is publicly available. See
            \cite{JonesData} and \cite{KhovanovData}.
        }
        This allows us to measure how common it is for a knot
        to have a unique Jones or Khovanov polynomial when compared with other
        low crossing knots. Unsurprisingly, the Jones polynomial does not
        distinguish nearly as many knots as the Khovanov polynomial.%
        \footnote{%
            The Khovanov polynomial throws away the torsion component of
            Khovanov homology. It would be interesting to compare how much
            stronger the \textit{full} Khovanov homology is.
        }
        The \textit{statistics} for these knot invariants are shown in
        Tabs.~\ref{tab:jones_stat} and \ref{tab:kho_stat}, the key for the
        tables is given in Tab.~\ref{tab:key}.
        \par\hfill\par
        Computational limits (temporarily) prohibit us from expanding
        Tab.~\ref{tab:kho_stat} to 20 crossings, but this is being worked on.
        After examining the data one sees that the Khovanov polynomial
        distinguishes about $32\%$ of knots with up to 17 crossings. Put it
        another way, about 7 out of 10 knots share their Khovanov polynomial
        with another knot in the list.
        \par\hfill\par
        We can use this data to play around with the statistics and ask whether or
        not the conjectures we've made in the previous sections are sound. That is,
        we can now provide a heuristic argument. The following arguments were
        greatly aided by collaborations between the second author and
        Peter Doyle.
        \par\hfill\par
        The smaller the probability of a knot having a unique Khovanov
        polynomial among all other prime knots up to $N$ crossings, the
        more likely our findings are statistically significant. We note that
        the probability of having a unique Jones or Khovanov polynomial appears
        to be monotonically \textit{decreasing} with crossing number. Since
        we do not have the Khovanov data for 18 to 20 crossings, we use
        the value for 17 crossings. If the trend continues and the probability
        for $N=20$ is smaller, that will only make our heuristic argument
        stronger. To avoid mixing the statistics, we also use the 17 crossing
        probability for the Jones polynomial. Here we \textit{know} that
        20 crossings produces a smaller probability since this data is
        available to us.
        \par\hfill\par
        A first approach is to say that the probability of a prime knot having a
        unique Khovanov polynomial among all other prime
        knots with up to 17 crossings is 0.32 (Tab.~\ref{tab:kho_stat}).
        Ignoring mirrors, there are 16 twist knots and 14 torus knots with
        less than or equal to 17 crossings. For twist knots these are the knots
        $K_{m}$ for $m=0,\,1,\,\cdots,\,15$, and for torus knots these are given
        by the co-prime pairs
        $(2, 3)$, $(2, 5)$, $(2, 7)$, $(2, 9)$, $(2, 11)$, $(2, 13)$,
        $(2, 15)$, $(2, 17)$, $(3, 4)$, $(3, 5)$, $(3, 7)$, $(3, 8)$, $(4, 5)$,
        and the unknot.
        The (na\"{i}ve) probability of all of these knots having unique
        Khovanov polynomials is $0.32^{16+14-1}\approx{4.5}\times{10}^{-15}$, where
        the $-1$ in the exponent comes from double counting the unknot, once as
        a twist knot and once as a torus knot.
        For twist knots alone this is $0.32^{16}\approx{1.2}\times{10}^{-8}$ and
        for torus knots this is $0.32^{14}\approx{1.2}\times{10}^{-7}$.
        \par\hfill\par
        These are ridiculously small probabilities, which is good for our
        conjecture,
        but this assumes these knots are more-or-less \textit{random}.
        Let us use their Jones polynomials to demonstrate that they more-or-less
        are. The Jones polynomial
        has a $0.26$ probability of distinguishing a knot in our list of prime knots
        up to 17 crossings. Since having a unique Jones polynomial would imply
        having a unique Khovanov polynomial (recall $Kh_{K}(q,\,-1)=J_{K}(q)$),
        the set of knots
        that are distinguished by their Jones polynomial is a subset of the set of
        knots distinguished by the Khovanov polynomial. If we label
        \texttt{KH} as the event that a knot is determined by its Khovanov
        polynomial (when considered among prime knots with 17 crossings or less),
        and similarly label \texttt{J} for the Jones polynomial, the conditional
        probability is:
        \begin{equation}
            P(\textrm{not }\texttt{J}\;|\;\texttt{KH})
            =\frac{P(\texttt{KH})-P(\texttt{J})}{P(\texttt{KH})}
            \approx{0.21}
        \end{equation}
        So with a $0.21$ probability a knot that is known to be distinguished by
        its Khovanov polynomial will not be distinguished by its Jones polynomial.
        We have shown by brute force methods that the torus and twist knots are
        part of the Khovanov-distinguished family for prime knots up to 17
        crossings, so we can apply this probability. There were 4 torus knots that
        are not distinguished by their Jones polynomial and 9 such twist knots.%
        \footnote{%
            Tab.~\ref{table:matching_twist_knots} has 11 such twist knots, but
            two of these matches are for knots with more than 17 crossings. To
            avoid using the wrong probabilities with this data, we stick to
            17 crossings.
        }
        When accounting for mirrors, this gives a
        $0.15$ probability of a torus knot being distinguished by the Khovanov
        polynomial and not the Jones polynomial, and a $0.3$ probability for twist
        knots. Roughly, this means torus knots are more likely to be distinguished
        by the Jones polynomial than the average Khovanov-distinguished knot,
        and the twist knots are less likely. The total probability for the union
        of these two families is $0.23$, which is close to the general
        probability of $0.21$. It is then fair to say that these
        families are behaving \textit{randomly}, as far as their
        Jones detection rates are concerned.
        \par\hfill\par
        We can also look at the conditional probability of finding a knot that is
        distinguished by its Khovanov polynomial given that we know it is not
        distinguished by its Jones polynomial. First, a quick review of probability.
        Given events $A,B\subseteq{X}$, for some probability space $X$, if $P(B)$
        is non-zero, then:
        \begin{align}
            P(B)&=P\big((A\cap{B})\cup(A^{C}\cap{B})\big)\\
            &=P(A\cap{B})+P(A^{C}\cap{B})\\
            \Rightarrow
            1&=\frac{P(A\cap{B})}{P(B)}+\frac{P(A^{C}\cap{B})}{P(B)}\\
            &=P(A\;|\;B)+P(A^{C}\;|\;B)
        \end{align}
        where $A^{C}=X\setminus{A}$, the complement of $A$ in $X$, and
        where $P(A\;|\;B)$ denotes the conditional probability. Note that
        $P(A^{C})$ is the same thing as $P(\textrm{not }A)$.
        Rearranging we have:
        \begin{equation}
            P(A\;|\;B)=1-P(A^{C}\;|\;B)
        \end{equation}
        a rather standard and intuitive result, we will now use this in our
        following argument. We seek the probability
        $P(\texttt{KH}\;|\;\textrm{not }\texttt{J})$.
        That is, the probability a knot will
        be distinguished by its Khovanov polynomial given that we know it is
        not distinguished by its Jones polynomial. Combining our previous derivation
        with Bayes' theorem, we obtain:
        \begin{align}
            P(\texttt{KH}\;|\;\textrm{not }\texttt{J})
            &=P(\textrm{not }\texttt{J}\;|\;\texttt{KH})
            \frac{P(\texttt{KH})}{P(\textrm{not }\texttt{J})}\\
            &=\big(1-P(\texttt{J}\;|\;\texttt{KH})\big)
                \frac{P(\texttt{KH})}{P(\textrm{not }\texttt{J})}\\
            &=\big(
                1-\frac{P(\texttt{KH}\;|\;\textrm{J})P(\texttt{J})}{P(\texttt{KH})}
            \big)
            \frac{P(\texttt{KH})}{P(\textrm{not }\texttt{J})}\\
            &=\big(
                1-\frac{P(\texttt{J})}{P(\texttt{KH})}
            \big)
            \frac{P(\texttt{KH})}{P(\textrm{not }\texttt{J})}\\
            &=\big(
                \frac{P(\texttt{KH})-P(\texttt{J})}{P(\texttt{KH})}
            \big)
            \frac{P(\texttt{KH})}{P(\textrm{not }\texttt{J})}\\
            &=\frac{P(\texttt{KH})-P(\texttt{J})}{P(\textrm{not }\texttt{J})}\\
            &=\frac{P(\texttt{KH})-P(\texttt{J})}{1-P(\texttt{J})}\\
            &=0.092
        \end{align}
        We have used the fact that $P(\texttt{KH}\;|\;\texttt{J})=1$ to obtain our
        result.
        To be on the conservative side we can round this up (smaller probabilities
        make our conjecture stronger) to $10\%$. So with roughly a $10\%$
        probability a knot that is not Jones-distinguished will still be
        distinguished by its Khovanov polynomial. We found 11 twist knots and 4
        torus knots that are not Jones-distinguished, but are still
        Khovanov-distinguished. Hence for the torus knots we may choose a $p$
        value of $0.001$, and for the twist knots we can go much lower, indicating
        that it is likely a statistically significant result.
        \par\hfill\par
        Now one may say that \textit{most} of the prime knots with 17 or less
        crossings happen to have 17 crossings, whereas most of the twist knots
        we've examined do not.
        A similar claim can be made for the torus knots. We can reformulate
        our probability argument as follows to account for this. Let
        $p(k)$ be the fraction of prime knots of up to $k$ crossings that are
        distinguished by their Khovanov polynomial. Let
        $Cr(K)$ denote the crossing number of $K$. The product
        $\prod_{m}p\big(Cr(K_{m})\big)$ over all twist knots with up to 17 crossings
        then serves as a better probability that all should be
        Khovanov-distinguished. Using our tables this number is
        $2.8\times{10}^{-5}$, which is still very small. For the torus
        knots we obtain $8.4\times{10}^{-6}$, a tiny probability.
        \par\hfill\par
        Lastly, one may argue that several of the torus and twist knots are known
        to be Khovanov-distinguished among all possible knots, like the unknot,
        and so should not be included in these probabilities.
        Our previous argument handles this since these knots have
        crossing number no greater than 5, and $p(k)=1$ for $k\leq{9}$.
        \par\hfill\par
        The premise behind these arguments, that knots are random, is a point of
        contention. Nevertheless the above argument is not entirely meaningless
        and the authors believe it does motivate further study of these conjectures.
        \begin{table}
            \centering
            \begin{tabular}{| l | l |}
                \hline
                Keyword & Description\\
                \hline
                Cr     & Crossing number, largest number of crossings considered.\\
                Unique & Number of polynomials that occur for one knot.\\
                Almost & Number of polynomials that occur for two knots.\\
                Total  & Total number of distinct polynomials in list.\\
                Knots  & Total number of knots in list.\\
                FracU  & Unique / Total\\
                FracT  & Total / Knots\\
                FracK  & Unique / Knots\\
                \hline
            \end{tabular}
            \caption{Legend for Tabs.~\ref{tab:jones_stat} and \ref{tab:kho_stat}}
            \label{tab:key}
        \end{table}
        \begin{table}
            \centering
            \begin{tabular}{| r | r | r | r | r | r | r | r |}
                \hline
                Cr &  Unique  &  Almost  &   Total   &   Knots    &  FracU   &  FracT   &  FracK\\
                \hline
                03 &        1 &        0 &         1 &         1 & 1.000000 & 1.000000 & 1.000000\\
                04 &        2 &        0 &         2 &         2 & 1.000000 & 1.000000 & 1.000000\\
                05 &        4 &        0 &         4 &         4 & 1.000000 & 1.000000 & 1.000000\\
                06 &        7 &        0 &         7 &         7 & 1.000000 & 1.000000 & 1.000000\\
                07 &       14 &        0 &        14 &        14 & 1.000000 & 1.000000 & 1.000000\\
                08 &       35 &        0 &        35 &        35 & 1.000000 & 1.000000 & 1.000000\\
                09 &       84 &        0 &        84 &        84 & 1.000000 & 1.000000 & 1.000000\\
                10 &      223 &       13 &       236 &       249 & 0.944915 & 0.947791 & 0.895582\\
                11 &      626 &       77 &       710 &       801 & 0.881690 & 0.886392 & 0.781523\\
                12 &     1981 &      345 &      2420 &      2977 & 0.818595 & 0.812899 & 0.665435\\
                13 &     6855 &     1695 &      9287 &     12965 & 0.738129 & 0.716313 & 0.528731\\
                14 &    25271 &     7439 &     37578 &     59937 & 0.672495 & 0.626958 & 0.421626\\
                15 &   105246 &    35371 &    170363 &    313230 & 0.617775 & 0.543891 & 0.336002\\
                16 &   487774 &   173677 &    829284 &   1701935 & 0.588187 & 0.487260 & 0.286600\\
                17 &  2498968 &   894450 &   4342890 &   9755328 & 0.575416 & 0.445181 & 0.256164\\
                18 & 13817237 &  4863074 &  24116048 &  58021794 & 0.572948 & 0.415638 & 0.238139\\
                19 & 82712788 & 27409120 & 141439472 & 352152252 & 0.584793 & 0.401643 & 0.234878\\
                \hline
            \end{tabular}
            \caption{Statistics for the Jones Polynomial}
            \label{tab:jones_stat}
        \end{table}
        \begin{table}
            \centering
            \begin{tabular}{| r | r | r | r | r | r | r | r |}
                \hline
                Cr & Unique  & Almost  &  Total  &  Knots  &  FracU   &  FracT   &  FracK\\
                \hline
                03 &       1 &       0 &       1 &       1 & 1.000000 & 1.000000 & 1.000000\\
                04 &       2 &       0 &       2 &       2 & 1.000000 & 1.000000 & 1.000000\\
                05 &       4 &       0 &       4 &       4 & 1.000000 & 1.000000 & 1.000000\\
                06 &       7 &       0 &       7 &       7 & 1.000000 & 1.000000 & 1.000000\\
                07 &      14 &       0 &      14 &      14 & 1.000000 & 1.000000 & 1.000000\\
                08 &      35 &       0 &      35 &      35 & 1.000000 & 1.000000 & 1.000000\\
                09 &      84 &       0 &      84 &      84 & 1.000000 & 1.000000 & 1.000000\\
                10 &     225 &      12 &     237 &     249 & 0.949367 & 0.951807 & 0.903614\\
                11 &     641 &      71 &     718 &     801 & 0.892758 & 0.896380 & 0.800250\\
                12 &    2051 &     326 &    2462 &    2977 & 0.833063 & 0.827007 & 0.688949\\
                13 &    7223 &    1636 &    9539 &   12965 & 0.757207 & 0.735750 & 0.557115\\
                14 &   27317 &    7441 &   39222 &   59937 & 0.696471 & 0.654387 & 0.455762\\
                15 &  118534 &   36867 &  182598 &  313230 & 0.649153 & 0.582952 & 0.378425\\
                16 &  578928 &  187639 &  919835 & 1701935 & 0.629382 & 0.540464 & 0.340159\\
                17 & 3167028 & 1001101 & 5033403 & 9755328 & 0.629202 & 0.515965 & 0.324646\\
                \hline
            \end{tabular}
            \caption{Statistics for the Khovanov Polynomial}
            \label{tab:kho_stat}
        \end{table}
    \section{Acknowledgements}
        We thank Shadi Ali Ahmad, Ilya Kryukov and Jacob Swenberg for their
        observation and remarks in the early stage of this work. The numerical
        results of Samantha Allen and Jacob Swenberg about the Jones Polynomial
        of a connected sum of two Hopf links were very inspirational since they
        provided the connection to the work of Fan Ding and Hansjorg Geiges on
        the Legendrian simple links. We are very grateful to Nikolay Pultsin
        for his indispensable help with the Khovanov homology code available
        on the internet. We are thankful to Jim Hoste, Morwen Thistlethwaite
        and Jeff Weeks for sharing their database of DT codes of knots, and
        to Benjamin Burton for providing his 19 crossing database freely to
        the public. The algorithm that was devised to compute the Jones
        polynomial efficiently (the software is publicly available,
        see \cite{MaguireLibtmpl}) was made possible by contributions from
        Peter Doyle, and we also thank him for guiding us in our statistical
        argument. We thank Steven Sivek and Matthew Heddon
        for alarming us about the bug that existed in our KFH code.
    \newpage
    \section{Appendix A: Torus Knot Data}
        \begin{table}[H]
            \centering
            \begin{tabular}{| c | c | c | c | c | c | c | c |}
                \hline
                $q\symbol{92}t$&$-5$&$-4$&$-3$&$-2$&$-1$&$0$&$\chi$\\
                \hline
                $-15$&1&&&&&&$-1$\\
                \hline
                $-13$&&&&&&&\\
                \hline
                $-11$&&1&1&&&&\\
                \hline
                $-9$&&&&&&&\\
                \hline
                $-7$&&&&1&&&1\\
                \hline
                $-5$&&&&&&1&1\\
                \hline
                $-3$&&&&&&1&1\\
                \hline
            \end{tabular}
            \caption{Khovanov Polynomial for $T(5,2)$}
            \label{table:t_2_5_kho}
        \end{table}
        \begin{table}[H]
            \centering
            \begin{tabular}{| c | c | c | c | c | c | c | c | c | c |}
                \hline
                $q\symbol{92}t$&$-7$&$-6$&$-5$&$-4$&$-3$&$-2$&$-1$&$0$&$\chi$\\
                \hline
                $-15$&1&&&&&&&&$-1$\\
                \hline
                $-13$&&&&&&&&&\\
                \hline
                $-11$&&1&1&&&&&&\\
                \hline
                $-9$&&&&1&1&&&&\\
                \hline
                $-7$&&&&1&&&&&1\\
                \hline
                $-5$&&&&&1&2&&&1\\
                \hline
                $-3$&&&&&&&&1&1\\
                \hline
                $-1$&&&&&&&1&1&\\
                \hline
            \end{tabular}
            \caption{Khovanov Polynomial for \texttt{dciaFHjEbg}}
            \label{table:dciaFHjEbg_kho}
        \end{table}
        \begin{table}[H]
            \centering
            \begin{tabular}{| c | c | c | c | c | c | c | c | c | c | c | c |}
                \hline
                $q\symbol{92}t$&$-9$&$-8$&$-7$&$-6$&$-5$&$-4$&$-3$&$-2$&$-1$&$0$&$\chi$\\
                \hline
                $-15$&1&&&&&&&&&&$-1$\\
                \hline
                $-13$&&&&&&&&&&&\\
                \hline
                $-11$&&1&1&&&&&&&&\\
                \hline
                $-9$&&&&1&1&&&&&&\\
                \hline
                $-7$&&&&1&&&&&&&1\\
                \hline
                $-5$&&&&&1&2&&&&&1\\
                \hline
                $-3$&&&&&&&&1&&&1\\
                \hline
                $-1$&&&&&&&1&&&1&\\
                \hline
                $1$&&&&&&&&&1&1&\\
                \hline
            \end{tabular}
            \caption{Khovanov Polynomial for \texttt{iNHlPJqCoKFmdABgE}}
            \label{table:iNHlPJqCoKFmdABgE_kho}
        \end{table}
        \begin{table}[H]
            \centering
            \begin{tabular}{| c | c | c | c | c | c | c | c | c | c |}
                \hline
                $q\symbol{92}t$&$-7$&$-6$&$-5$&$-4$&$-3$&$-2$&$-1$&$0$&$\chi$\\
                \hline
                $-21$&1&&&&&&&&$-1$\\
                \hline
                $-19$&&&&&&&&&\\
                \hline
                $-17$&&1&1&&&&&&\\
                \hline
                $-15$&&&&&&&&&\\
                \hline
                $-13$&&&&1&1&&&&\\
                \hline
                $-11$&&&&&&&&&\\
                \hline
                $-9$&&&&&&1&&&1\\
                \hline
                $-7$&&&&&&&&1&1\\
                \hline
                $-5$&&&&&&&&1&1\\
                \hline
            \end{tabular}
            \caption{Khovanov Polynomial for $T(7,2)$}
            \label{table:t_7_2_kho}
        \end{table}
        \begin{table}[H]
            \centering
            \begin{tabular}{| c | c | c | c | c | c | c | c | c | c | c | c |}
                \hline
                $q\symbol{92}t$&$-9$&$-8$&$-7$&$-6$&$-5$&$-4$&$-3$&$-2$&$-1$&$0$&$\chi$\\
                \hline
                $-21$&1&&&&&&&&&&$-1$\\
                \hline
                $-19$&&&&&&&&&&&\\
                \hline
                $-17$&&1&1&&&&&&&&\\
                \hline
                $-15$&&&&1&1&&&&&&\\
                \hline
                $-13$&&&&1&1&&&&&&\\
                \hline
                $-11$&&&&&1&2&1&&&&\\
                \hline
                $-9$&&&&&&1&&&&&1\\
                \hline
                $-7$&&&&&&&1&2&&&1\\
                \hline
                $-5$&&&&&&&&&&1&1\\
                \hline
                $-3$&&&&&&&&&1&1&\\
                \hline
            \end{tabular}
            \caption{Khovanov Polynomial for \texttt{fJGkHlICEABd}}
            \label{table:fJGkHlICEABd_kho}
        \end{table}
        \begin{table}[H]
            \centering
            \begin{tabular}{| c | c | c | c | c | c | c | c | c | c | c | c | c | c |}
                \hline
                $q\symbol{92}t$&$0$&$1$&$2$&$3$&$4$&$5$&$6$&$7$&$8$&$9$&$10$&$11$&$\chi$\\
                \hline
                $9$&1&&&&&&&&&&&&1\\
                \hline
                $11$&1&&&&&&&&&&&&1\\
                \hline
                $13$&&&1&&&&&&&&&&1\\
                \hline
                $15$&&&&&&&&&&&&&\\
                \hline
                $17$&&&&1&1&&&&&&&&\\
                \hline
                $19$&&&&&&&&&&&&&\\
                \hline
                $21$&&&&&&1&1&&&&&&\\
                \hline
                $23$&&&&&&&&&&&&&\\
                \hline
                $25$&&&&&&&&1&1&&&&\\
                \hline
                $27$&&&&&&&&&&&&&\\
                \hline
                $29$&&&&&&&&&&1&1&&\\
                \hline
                $31$&&&&&&&&&&&&&\\
                \hline
                $33$&&&&&&&&&&&&1&1\\
                \hline
            \end{tabular}
            \caption{Khovanov Polynomial for $T(11,\,2)$}
            \label{table:t_2_11_kho}
        \end{table}
        \begin{table}[H]
            \centering
            \begin{tabular}{| c | c | c | c | c | c | c | c | c | c | c | c | c | c | c | c |}
                \hline
                $q\symbol{92}t$&$0$&$1$&$2$&$3$&$4$&$5$&$6$&$7$&$8$&$9$&$10$&$11$&$12$&$13$&$\chi$\\
                \hline
                $9$&1&&&&&&&&&&&&&&1\\
                \hline
                $11$&1&&&&&&&&&&&&&&1\\
                \hline
                $13$&&&1&&&&&&&&&&&&1\\
                \hline
                $15$&&&&&1&1&&&&&&&&&\\
                \hline
                $17$&&&&1&1&&&&&&&&&&\\
                \hline
                $19$&&&&&&1&2&1&&&&&&&\\
                \hline
                $21$&&&&&&1&&&1&&&&&&\\
                \hline
                $23$&&&&&&&&1&2&1&&&&&\\
                \hline
                $25$&&&&&&&&&&1&1&&&&\\
                \hline
                $27$&&&&&&&&&&1&1&&&&\\
                \hline
                $29$&&&&&&&&&&&&1&1&&\\
                \hline
                $31$&&&&&&&&&&&&&&&\\
                \hline
                $33$&&&&&&&&&&&&&&1&1\\
                \hline
            \end{tabular}
            \caption{Khovanov Polynomial for \texttt{gHlImJnKBDFAce}}
            \label{table:t_gHlImJnKBDFAce_kho}
        \end{table}
    \newpage
    \section{Appendix B: Twist Knot Data}
        \begin{table}[H]
            \centering
            \begin{tabular}{| c | c | c | c | c | c | c |}
                \hline
                $q\symbol{92}t$&$-2$&$-1$&$0$&$1$&$2$&$\chi$\\
                \hline
                $-5$&1&&&&&1\\
                \hline
                $-3$&&&&&&\\
                \hline
                $-1$&&1&1&&&\\
                \hline
                $1$&&&1&1&&\\
                \hline
                $3$&&&&&&\\
                \hline
                $5$&&&&&1&1\\
                \hline
            \end{tabular}
            \caption{Khovanov Polynomial for the Figure-Eight Knot}
            \label{table:m_2_kho}
        \end{table}
        \begin{table}[H]
            \centering
            \begin{tabular}{| c | c | c | c | c | c | c | c | c | c |}
                \hline
                $q\symbol{92}t$&$-4$&$-3$&$-2$&$-1$&$0$&$1$&$2$&$3$&$\chi$\\
                \hline
                $-5$&1&&&&&&&&1\\
                \hline
                $-3$&&&&&&&&&\\
                \hline
                $-1$&&1&1&&&&&&\\
                \hline
                $1$&&&&1&1&&&&\\
                \hline
                $3$&&&&1&1&&&&\\
                \hline
                $5$&&&&&1&1&1&&1\\
                \hline
                $7$&&&&&&&&&\\
                \hline
                $9$&&&&&&&1&1&\\
                \hline
            \end{tabular}
            \caption{Khovanov Polynomial for \texttt{eikGbHJCaFd}}
            \label{table:eikGbHJCaFd_kho}
        \end{table}
        \begin{table}[H]
            \centering
            \begin{tabular}{| c | c | c | c | c | c | c | c |}
                \hline
                $q\symbol{92}t$&$0$&$1$&$2$&$3$&$4$&$5$&$\chi$\\
                \hline
                $1$&1&&&&&&1\\
                \hline
                $3$&1&1&&&&&\\
                \hline
                $5$&&&1&&&&1\\
                \hline
                $7$&&&1&&&&1\\
                \hline
                $9$&&&&1&1&&\\
                \hline
                $11$&&&&&&&\\
                \hline
                $13$&&&&&&1&$-1$\\
                \hline
            \end{tabular}
            \caption{Khovanov Polynomial for the $5_{2}$ Knot}
            \label{table:m_3_kho}
        \end{table}
        \begin{table}[H]
            \centering
            \begin{tabular}{| c | c | c | c | c | c | c | c | c | c | c |}
                \hline
                $q\symbol{92}t$&$-2$&$-1$&$0$&$1$&$2$&$3$&$4$&$5$&$6$&$\chi$\\
                \hline
                $1$&1&&&&&&&&&1\\
                \hline
                $3$&&&&&&&&&&\\
                \hline
                $5$&&1&2&&&&&&&1\\
                \hline
                $7$&&&1&1&1&&&&&1\\
                \hline
                $9$&&&&1&1&&&&&\\
                \hline
                $11$&&&&&1&2&1&&&\\
                \hline
                $13$&&&&&&1&1&1&&$-1$\\
                \hline
                $15$&&&&&&&1&1&&\\
                \hline
                $17$&&&&&&&&1&1&\\
                \hline
            \end{tabular}
            \caption{Khovanov Polynomial for \texttt{dgikFHaEjbc}}
            \label{table:dgikFHaEjbc_kho}
        \end{table}
        \begin{table}[H]
            \centering
            \begin{tabular}{| c | c | c | c | c | c | c | c | c | c |}
                \hline
                $q\symbol{92}t$&$0$&$1$&$2$&$3$&$4$&$5$&$6$&$7$&$\chi$\\
                \hline
                $-1$&1&1&&&&&&&\\
                \hline
                $1$&1&&&&&&&&1\\
                \hline
                $3$&&&1&1&&&&&\\
                \hline
                $5$&&&&&1&&&&1\\
                \hline
                $7$&&&&&1&&&&1\\
                \hline
                $9$&&&&&&1&1&&\\
                \hline
                $11$&&&&&&&&&\\
                \hline
                $13$&&&&&&&&1&$-1$\\
                \hline
            \end{tabular}
            \caption{Khovanov Polynomial for \texttt{gfJKHlaIEBCD}}
            \label{table:gfJKHlaIEBCD_kho}
        \end{table}
        \begin{table}[H]
            \centering
            \begin{tabular}{| c | c | c | c | c | c | c | c | c | c | c | c | c |}
                \hline
                $q\symbol{92}t$&$-2$&$-1$&$0$&$1$&$2$&$3$&$4$&$5$&$6$&$7$&$8$&$\chi$\\
                \hline
                $1$&1&&&&&&&&&&&1\\
                \hline
                $3$&&&&&&&&&&&&\\
                \hline
                $5$&&1&2&&&&&&&&&1\\
                \hline
                $7$&&&1&1&1&&&&&&&1\\
                \hline
                $9$&&&&1&1&&&&&&&\\
                \hline
                $11$&&&&&1&2&1&&&&&\\
                \hline
                $13$&&&&&&1&1&1&&&&$-1$\\
                \hline
                $15$&&&&&&&1&1&&&&\\
                \hline
                $17$&&&&&&&&1&2&1&&\\
                \hline
                $19$&&&&&&&&&&&&\\
                \hline
                $21$&&&&&&&&&&1&1&\\
                \hline
            \end{tabular}
            \caption{Khovanov Polynomial for \texttt{hGJaMlCdEKBfI}}
            \label{table:hGJaMlCdEKBfI_kho}
        \end{table}
        \begin{table}[H]
            \centering
            \begin{tabular}{| c | c | c | c | c | c | c | c | c | c |}
                \hline
                $q\symbol{92}t$&$0$&$1$&$2$&$3$&$4$&$5$&$6$&$7$&$\chi$\\
                \hline
                $1$&1&&&&&&&&1\\
                \hline
                $3$&1&1&&&&&&&\\
                \hline
                $5$&&&1&&&&&&1\\
                \hline
                $7$&&&1&1&&&&&\\
                \hline
                $9$&&&&1&1&&&&\\
                \hline
                $11$&&&&&1&&&&1\\
                \hline
                $13$&&&&&&1&1&&\\
                \hline
                $15$&&&&&&&&&\\
                \hline
                $17$&&&&&&&&1&$-1$\\
                \hline
            \end{tabular}
            \caption{Khovanov Polynomial for the $7_{2}$ Knot}
            \label{table:m_5_kho}
        \end{table}
        \begin{table}[H]
            \centering
            \begin{tabular}{| c | c | c | c | c | c | c | c | c | c |}
                \hline
                $q\symbol{92}t$&$-2$&$-1$&$0$&$1$&$2$&$3$&$4$&$5$&$\chi$\\
                \hline
                $1$&1&&&&&&&&1\\
                \hline
                $3$&&&&&&&&&\\
                \hline
                $5$&&1&2&&&&&&1\\
                \hline
                $7$&&&1&1&&&&&\\
                \hline
                $9$&&&&1&1&&&&\\
                \hline
                $11$&&&&&1&1&1&&1\\
                \hline
                $13$&&&&&&1&1&&\\
                \hline
                $15$&&&&&&&1&1&\\
                \hline
                $17$&&&&&&&&1&$-1$\\
                \hline
            \end{tabular}
            \caption{Khovanov Polynomial for \texttt{bhDGijCkaef}}
            \label{table:bhDGijCkaef_kho}
        \end{table}
        \begin{table}[H]
            \centering
            \begin{tabular}{| c | c | c | c | c | c | c | c | c | c | c |}
                \hline
                $q\symbol{92}t$&$-2$&$-1$&$0$&$1$&$2$&$3$&$4$&$5$&$6$&$\chi$\\
                \hline
                $-5$&1&&&&&&&&&1\\
                \hline
                $-3$&&&&&&&&&&\\
                \hline
                $-1$&&1&2&&&&&&&1\\
                \hline
                $1$&&&1&1&&&&&&\\
                \hline
                $3$&&&&1&1&&&&&\\
                \hline
                $5$&&&&&1&1&&&&\\
                \hline
                $7$&&&&&&1&&&&$-1$\\
                \hline
                $9$&&&&&&&1&1&&\\
                \hline
                $11$&&&&&&&&&&\\
                \hline
                $13$&&&&&&&&&1&1\\
                \hline
            \end{tabular}
            \caption{Khovanov Polynomial for the $8_{1}$ Knot}
            \label{table:m_6_kho}
        \end{table}
        \begin{table}[H]
            \centering
            \begin{tabular}{| c | c | c | c | c | c | c | c | c | c | c |}
                \hline
                $q\symbol{92}t$&$-4$&$-3$&$-2$&$-1$&$0$&$1$&$2$&$3$&$4$&$\chi$\\
                \hline
                $-5$&1&&&&&&&&&1\\
                \hline
                $-3$&&&&&&&&&&\\
                \hline
                $-1$&&1&2&&&&&&&1\\
                \hline
                $1$&&&&1&1&&&&&\\
                \hline
                $3$&&&&2&2&&&&&\\
                \hline
                $5$&&&&&1&2&1&&&\\
                \hline
                $7$&&&&&&1&&&&$-1$\\
                \hline
                $9$&&&&&&&2&2&&\\
                \hline
                $11$&&&&&&&&&&\\
                \hline
                $13$&&&&&&&&&1&1\\
                \hline
            \end{tabular}
            \caption{Khovanov Polynomial for \texttt{cefIgbajkDh}}
            \label{table:cefIgbajkDh_kho}
        \end{table}
        \begin{table}[H]
            \centering
            \begin{tabular}{| c | c | c | c | c | c | c | c | c | c | c |}
                \hline
                $q\symbol{92}t$&$-4$&$-3$&$-2$&$-1$&$0$&$1$&$2$&$3$&$4$&$\chi$\\
                \hline
                $-5$&1&&&&&&&&&1\\
                \hline
                $-3$&&&&&&&&&&\\
                \hline
                $-1$&&1&2&&&&&&&1\\
                \hline
                $1$&&&&1&1&&&&&\\
                \hline
                $3$&&&&2&2&&&&&\\
                \hline
                $5$&&&&&1&2&1&&&\\
                \hline
                $7$&&&&&&1&&&&$-1$\\
                \hline
                $9$&&&&&&&2&2&&\\
                \hline
                $11$&&&&&&&&&&\\
                \hline
                $13$&&&&&&&&&1&1\\
                \hline
            \end{tabular}
            \caption{Khovanov Polynomial for \texttt{femIbaJKLCGHd}}
            \label{table:femIbaJKLCGHd_kho}
        \end{table}
        \begin{table}[H]
            \centering
            \begin{tabular}{| c | c | c | c | c | c | c | c | c | c | c | c | c | c | c |}
                \hline
                $q\symbol{92}t$&$-4$&$-3$&$-2$&$-1$&$0$&$1$&$2$&$3$&$4$&$5$&$6$&$7$&$8$&$\chi$\\
                \hline
                $-5$&1&&&&&&&&&&&&&1\\
                \hline
                $-3$&&&&&&&&&&&&&&\\
                \hline
                $-1$&&1&2&&&&&&&&&&&1\\
                \hline
                $1$&&&&1&1&&&&&&&&&\\
                \hline
                $3$&&&&2&2&&&&&&&&&\\
                \hline
                $5$&&&&&2&4&2&&&&&&&\\
                \hline
                $7$&&&&&&1&1&1&&&&&&$-1$\\
                \hline
                $9$&&&&&&&3&4&1&&&&&\\
                \hline
                $11$&&&&&&&&1&3&2&&&&\\
                \hline
                $13$&&&&&&&&&2&1&&&&1\\
                \hline
                $15$&&&&&&&&&&2&3&1&&\\
                \hline
                $17$&&&&&&&&&&&&&&\\
                \hline
                $19$&&&&&&&&&&&&1&1&\\
                \hline
            \end{tabular}
            \caption{Khovanov Polynomial for \texttt{jpIFNMrClqOhkEDabg}}
            \label{table:jpIFNMrClqOhkEDabg_kho}
        \end{table}
        \begin{table}[H]
            \centering
            \begin{tabular}{| c | c | c | c | c | c | c | c | c | c | c | c |}
                \hline
                $q\symbol{92}t$&$0$&$1$&$2$&$3$&$4$&$5$&$6$&$7$&$8$&$9$&$\chi$\\
                \hline
                $1$&1&&&&&&&&&&1\\
                \hline
                $3$&1&1&&&&&&&&&\\
                \hline
                $5$&&&1&&&&&&&&1\\
                \hline
                $7$&&&1&1&&&&&&&\\
                \hline
                $9$&&&&1&1&&&&&&\\
                \hline
                $11$&&&&&1&1&&&&&\\
                \hline
                $13$&&&&&&1&1&&&&\\
                \hline
                $15$&&&&&&&1&&&&1\\
                \hline
                $17$&&&&&&&&1&1&&\\
                \hline
                $19$&&&&&&&&&&&\\
                \hline
                $21$&&&&&&&&&&1&$-1$\\
                \hline
            \end{tabular}
            \caption{Khovanov Polynomial for $K_{7}$}
            \label{table:m_7_kho}
        \end{table}
        \begin{table}[H]
            \centering
            \begin{tabular}{| c | c | c | c | c | c | c | c | c | c | c | c |}
                \hline
                $q\symbol{92}t$&$-2$&$-1$&$0$&$1$&$2$&$3$&$4$&$5$&$6$&$7$&$\chi$\\
                \hline
                $1$&1&&&&&&&&&&1\\
                \hline
                $3$&&&&&&&&&&&\\
                \hline
                $5$&&1&2&&&&&&&&1\\
                \hline
                $7$&&&1&1&&&&&&&\\
                \hline
                $9$&&&&1&1&&&&&&\\
                \hline
                $11$&&&&&1&1&&&&&\\
                \hline
                $13$&&&&&&1&1&&&&\\
                \hline
                $15$&&&&&&&1&&&&1\\
                \hline
                $17$&&&&&&&&1&1&&\\
                \hline
                $19$&&&&&&&&&&&\\
                \hline
                $21$&&&&&&&&&&1&$-1$\\
                \hline
            \end{tabular}
            \caption{Khovanov Polynomial for \texttt{cgjFHIaDEkb}}
            \label{table:cgjFHIaDEkb_kho}
        \end{table}
        \begin{table}[H]
            \centering
            \begin{tabular}{| c | c | c | c | c | c | c | c | c | c | c | c | c |}
                \hline
                $q\symbol{92}t$&$-2$&$-1$&$0$&$1$&$2$&$3$&$4$&$5$&$6$&$7$&$8$&$\chi$\\
                \hline
                $-5$&1&&&&&&&&&&&1\\
                \hline
                $-3$&&&&&&&&&&&&\\
                \hline
                $-1$&&1&2&&&&&&&&&1\\
                \hline
                $1$&&&1&1&&&&&&&&\\
                \hline
                $3$&&&&1&1&&&&&&&\\
                \hline
                $5$&&&&&1&1&&&&&&\\
                \hline
                $7$&&&&&&1&1&&&&&\\
                \hline
                $9$&&&&&&&1&1&&&&\\
                \hline
                $11$&&&&&&&&1&&&&$-1$\\
                \hline
                $13$&&&&&&&&&1&1&&\\
                \hline
                $15$&&&&&&&&&&&&\\
                \hline
                $17$&&&&&&&&&&&1&1\\
                \hline
            \end{tabular}
            \caption{Khovanov Polynomial for $K_{8}$}
            \label{table:m_8_kho}
        \end{table}
        \begin{table}[H]
            \centering
            \begin{tabular}{| c | c | c | c | c | c | c | c | c | c | c | c | c | c | c |}
                \hline
                $q\symbol{92}t$&$-6$&$-5$&$-4$&$-3$&$-2$&$-1$&$0$&$1$&$2$&$3$&$4$&$5$&$6$&$\chi$\\
                \hline
                $-5$&1&&&&&&&&&&&&&1\\
                \hline
                $-3$&&&&&&&&&&&&&&\\
                \hline
                $-1$&&1&2&&&&&&&&&&&1\\
                \hline
                $1$&&&&1&1&&&&&&&&&\\
                \hline
                $3$&&&&2&1&&1&&&&&&&\\
                \hline
                $5$&&&&&1&3&2&&&&&&&\\
                \hline
                $7$&&&&&&1&&1&2&&&&&\\
                \hline
                $9$&&&&&&&2&2&&&&&&\\
                \hline
                $11$&&&&&&&&&1&3&1&&&$-1$\\
                \hline
                $13$&&&&&&&&&1&&&1&&\\
                \hline
                $15$&&&&&&&&&&&1&1&&\\
                \hline
                $17$&&&&&&&&&&&&&1&1\\
                \hline
            \end{tabular}
            \caption{Khovanov Polynomial for \texttt{knIHoBjCDQrMPaeLgF}}
            \label{table:knIHoBjCDQrMPaeLgF_kho}
        \end{table}
        \begin{table}[H]
            \centering
            \begin{tabular}{| c | c | c | c | c | c | c | c | c | c | c | c | c | c |}
                \hline
                $q\symbol{92}t$&$0$&$1$&$2$&$3$&$4$&$5$&$6$&$7$&$8$&$9$&$10$&$11$&$\chi$\\
                \hline
                $1$&1&&&&&&&&&&&&1\\
                \hline
                $3$&1&1&&&&&&&&&&&\\
                \hline
                $5$&&&1&&&&&&&&&&1\\
                \hline
                $7$&&&1&1&&&&&&&&&\\
                \hline
                $9$&&&&1&1&&&&&&&&\\
                \hline
                $11$&&&&&1&1&&&&&&&\\
                \hline
                $13$&&&&&&1&1&&&&&&\\
                \hline
                $15$&&&&&&&1&1&&&&&\\
                \hline
                $17$&&&&&&&&1&1&&&&\\
                \hline
                $19$&&&&&&&&&1&&&&1\\
                \hline
                $21$&&&&&&&&&&1&1&&\\
                \hline
                $23$&&&&&&&&&&&&&\\
                \hline
                $25$&&&&&&&&&&&&1&$-1$\\
                \hline
            \end{tabular}
            \caption{Khovanov Polynomial for $K_{9}$}
            \label{table:m_9_kho}
        \end{table}
        \begin{table}[H]
            \centering
            \begin{tabular}{| c | c | c | c | c | c | c | c | c | c | c | c | c | c |}
                \hline
                $q\symbol{92}t$&$-2$&$-1$&$0$&$1$&$2$&$3$&$4$&$5$&$6$&$7$&$8$&$9$&$\chi$\\
                \hline
                $1$&1&&&&&&&&&&&&1\\
                \hline
                $3$&&&&&&&&&&&&&\\
                \hline
                $5$&&1&2&&&&&&&&&&1\\
                \hline
                $7$&&&1&2&1&&&&&&&&\\
                \hline
                $9$&&&&1&1&&&&&&&&\\
                \hline
                $11$&&&&&2&3&1&&&&&&\\
                \hline
                $13$&&&&&&1&2&1&&&&&\\
                \hline
                $15$&&&&&&&2&2&&&&&\\
                \hline
                $17$&&&&&&&&2&3&1&&&\\
                \hline
                $19$&&&&&&&&&1&&&&1\\
                \hline
                $21$&&&&&&&&&&2&2&&\\
                \hline
                $23$&&&&&&&&&&&&&\\
                \hline
                $25$&&&&&&&&&&&&1&$-1$\\
                \hline
            \end{tabular}
            \caption{Khovanov Polynomial for \texttt{jopIFMrDlqNhkEabcg}}
            \label{table:jopIFMrDlqNhkEabcg_kho}
        \end{table}
    \newpage
    \section{Appendix C: Conjectured Legendrian Simple Knot Data}
        \begin{table}
            \centering
            \begin{tabular}{| c | c | c | c | c | c | c | c | c |}
                \hline
                $q\symbol{92}t$&$-2$&$-1$&$0$&$1$&$2$&$3$&$4$&$\chi$\\
                \hline
                $-3$&1&&&&&&&1\\
                \hline
                $-1$&&&&&&&&\\
                \hline
                $1$&&1&2&&&&&1\\
                \hline
                $3$&&&1&1&&&&\\
                \hline
                $5$&&&&1&1&&&\\
                \hline
                $7$&&&&&1&1&&\\
                \hline
                $9$&&&&&&1&&$-1$\\
                \hline
                $11$&&&&&&&1&1\\
                \hline
            \end{tabular}
            \caption{Khovanov Polynomial for $6_{2}$}
        \end{table}
        \begin{table}
            \centering
            \begin{tabular}{| c | c | c | c | c | c | c | c | c | c | c | c |}
                \hline
                $q\symbol{92}t$&$-4$&$-3$&$-2$&$-1$&$0$&$1$&$2$&$3$&$4$&$5$&$\chi$\\
                \hline
                $-3$&1&&&&&&&&&&1\\
                \hline
                $-1$&&&&&&&&&&&\\
                \hline
                $1$&&1&2&&&&&&&&1\\
                \hline
                $3$&&&&1&1&&&&&&\\
                \hline
                $5$&&&&2&2&&&&&&\\
                \hline
                $7$&&&&&1&2&1&&&&\\
                \hline
                $9$&&&&&&1&&&&&$-1$\\
                \hline
                $11$&&&&&&&2&2&1&&1\\
                \hline
                $13$&&&&&&&&&&&\\
                \hline
                $15$&&&&&&&&&1&1&\\
                \hline
            \end{tabular}
            \caption{Khovanov Polynomial for \texttt{glfoJcbKMNDaHIe}}
        \end{table}
        \begin{table}
            \centering
            \begin{tabular}{| c | c | c | c | c | c | c | c | c | c | c | c |}
                \hline
                $q\symbol{92}t$&$-4$&$-3$&$-2$&$-1$&$0$&$1$&$2$&$3$&$4$&$5$&$\chi$\\
                \hline
                $-3$&1&&&&&&&&&&1\\
                \hline
                $-1$&&&&&&&&&&&\\
                \hline
                $1$&&1&2&&&&&&&&1\\
                \hline
                $3$&&&&1&1&&&&&&\\
                \hline
                $5$&&&&2&2&&&&&&\\
                \hline
                $7$&&&&&1&2&1&&&&\\
                \hline
                $9$&&&&&&1&&&&&$-1$\\
                \hline
                $11$&&&&&&&2&2&1&&1\\
                \hline
                $13$&&&&&&&&&&&\\
                \hline
                $15$&&&&&&&&&1&1&\\
                \hline
            \end{tabular}
            \caption{Khovanov Polynomial for \texttt{hknEGmDbJLaIfc}}
        \end{table}
        \begin{table}
            \centering
            \begin{tabular}{| c | c | c | c | c | c | c | c | c | c | c | c | c |}
                \hline
                $q\symbol{92}t$&$-4$&$-3$&$-2$&$-1$&$0$&$1$&$2$&$3$&$4$&$5$&$6$&$\chi$\\
                \hline
                $-9$&1&1&&&&&&&&&&\\
                \hline
                $-7$&&&&&&&&&&&&\\
                \hline
                $-5$&&1&2&1&&&&&&&&\\
                \hline
                $-3$&&&&1&2&&&&&&&1\\
                \hline
                $-1$&&&&1&2&1&&&&&&\\
                \hline
                $1$&&&&&1&2&2&&&&&\\
                \hline
                $3$&&&&&&&1&1&&&&\\
                \hline
                $5$&&&&&&&1&2&1&&&\\
                \hline
                $7$&&&&&&&&&1&1&&\\
                \hline
                $9$&&&&&&&&&&1&&$-1$\\
                \hline
                $11$&&&&&&&&&&&1&1\\
                \hline
            \end{tabular}
            \caption{Khovanov Polynomial for \texttt{gKHlmIdJCEABf}}
        \end{table}
        \begin{table}
            \centering
            \begin{tabular}{| c | c | c | c | c | c | c | c | c | c | c | c | c |}
                \hline
                $q\symbol{92}t$&$-4$&$-3$&$-2$&$-1$&$0$&$1$&$2$&$3$&$4$&$5$&$6$&$\chi$\\
                \hline
                $-3$&1&&&&&&&&&&&1\\
                \hline
                $-1$&&&&&&&&&&&&\\
                \hline
                $1$&&1&2&&&&&&&&&1\\
                \hline
                $3$&&&&1&1&&&&&&&\\
                \hline
                $5$&&&&2&2&&&&&&&\\
                \hline
                $7$&&&&&2&3&1&&&&&\\
                \hline
                $9$&&&&&&1&1&1&&&&$-1$\\
                \hline
                $11$&&&&&&&2&2&1&&&1\\
                \hline
                $13$&&&&&&&&1&2&1&&\\
                \hline
                $15$&&&&&&&&&1&1&&\\
                \hline
                $17$&&&&&&&&&&1&1&\\
                \hline
            \end{tabular}
            \caption{Khovanov Polynomial for \texttt{ehkjmGIaFlcbd}}
        \end{table}
        \begin{table}
            \centering
            \begin{tabular}{| c | c | c | c | c | c | c | c | c | c |}
                \hline
                $q\symbol{92}t$&$0$&$1$&$2$&$3$&$4$&$5$&$6$&$7$&$\chi$\\
                \hline
                $3$&1&&&&&&&&1\\
                \hline
                $5$&1&1&&&&&&&\\
                \hline
                $7$&&&1&&&&&&1\\
                \hline
                $9$&&&1&1&&&&&\\
                \hline
                $11$&&&&1&2&&&&1\\
                \hline
                $13$&&&&&1&&&&1\\
                \hline
                $15$&&&&&&2&1&&$-1$\\
                \hline
                $17$&&&&&&&&&\\
                \hline
                $19$&&&&&&&&1&$-1$\\
                \hline
            \end{tabular}
            \caption{Khovanov Polynomial for $7_{3}$}
        \end{table}
        \begin{table}
            \centering
            \begin{tabular}{| c | c | c | c | c | c | c | c | c | c | c | c |}
                \hline
                $q\symbol{92}t$&$0$&$1$&$2$&$3$&$4$&$5$&$6$&$7$&$8$&$9$&$\chi$\\
                \hline
                $1$&1&1&&&&&&&&&\\
                \hline
                $3$&1&&&&&&&&&&1\\
                \hline
                $5$&&&2&2&&&&&&&\\
                \hline
                $7$&&&&&1&&&&&&1\\
                \hline
                $9$&&&&1&2&1&&&&&\\
                \hline
                $11$&&&&&&1&2&&&&1\\
                \hline
                $13$&&&&&&&1&&&&1\\
                \hline
                $15$&&&&&&&&2&1&&$-1$\\
                \hline
                $17$&&&&&&&&&&&\\
                \hline
                $19$&&&&&&&&&&1&$-1$\\
                \hline
            \end{tabular}
            \caption{Khovanov Polynomial for \texttt{hgelkIbaJFcd}}
        \end{table}
        \begin{table}
            \centering
            \begin{tabular}{| c | c | c | c | c | c | c | c | c | c |}
                \hline
                $q\symbol{92}t$&$0$&$1$&$2$&$3$&$4$&$5$&$6$&$7$&$\chi$\\
                \hline
                $1$&1&&&&&&&&1\\
                \hline
                $3$&1&2&&&&&&&$-1$\\
                \hline
                $5$&&&1&&&&&&1\\
                \hline
                $7$&&&2&1&&&&&1\\
                \hline
                $9$&&&&1&2&&&&1\\
                \hline
                $11$&&&&&1&&&&1\\
                \hline
                $13$&&&&&&2&1&&$-1$\\
                \hline
                $15$&&&&&&&&&\\
                \hline
                $17$&&&&&&&&1&$-1$\\
                \hline
            \end{tabular}
            \caption{Khovanov Polynomial for $7_{4}$}
        \end{table}
        \begin{table}
            \centering
            \begin{tabular}{| c | c | c | c | c | c | c | c | c | c | c |}
                \hline
                $q\symbol{92}t$&$-2$&$-1$&$0$&$1$&$2$&$3$&$4$&$5$&$6$&$\chi$\\
                \hline
                $1$&1&&&&&&&&&1\\
                \hline
                $3$&&1&&&&&&&&$-1$\\
                \hline
                $5$&&1&2&&&&&&&1\\
                \hline
                $7$&&&2&2&1&&&&&1\\
                \hline
                $9$&&&&1&2&&&&&1\\
                \hline
                $11$&&&&&2&2&1&&&1\\
                \hline
                $13$&&&&&&2&2&1&&$-1$\\
                \hline
                $15$&&&&&&&1&1&&\\
                \hline
                $17$&&&&&&&&2&1&1\\
                \hline
            \end{tabular}
            \caption{Khovanov Polynomial for \texttt{gfkHlbjIDaec}}
        \end{table}
        \begin{table}
            \centering
            \begin{tabular}{| c | c | c | c | c | c | c | c | c | c |}
                \hline
                $q\symbol{92}t$&$-5$&$-4$&$-3$&$-2$&$-1$&$0$&$1$&$2$&$\chi$\\
                \hline
                $-13$&2&&&&&&&&$-2$\\
                \hline
                $-11$&&1&&&&&&&1\\
                \hline
                $-9$&&2&3&&&&&&$-1$\\
                \hline
                $-7$&&&1&3&&&&&2\\
                \hline
                $-5$&&&&3&1&&&&2\\
                \hline
                $-3$&&&&&3&3&&&\\
                \hline
                $-1$&&&&&&2&1&&1\\
                \hline
                $1$&&&&&&&2&&$-2$\\
                \hline
                $3$&&&&&&&&1&1\\
                \hline
            \end{tabular}
            \caption{Khovanov Polynomial for $9_{48}$}
        \end{table}
        \begin{table}
            \centering
            \begin{tabular}{| c | c |c | c | c | c | c | c | c | c | c | c | c | c |}
                \hline
                $q\symbol{92}t$&$-7$&$-6$&$-5$&$-4$&$-3$&$-2$&$-1$&$0$&$1$&$2$&$3$&$4$&$\chi$\\
                \hline
                $-15$&1&1&&&&&&&&&&&\\
                \hline
                $-13$&&&2&&&&&&&&&&$-2$\\
                \hline
                $-11$&&1&2&2&&&&&&&&&1\\
                \hline
                $-9$&&&&3&5&1&&&&&&&$-1$\\
                \hline
                $-7$&&&&1&2&3&&&&&&&2\\
                \hline
                $-5$&&&&&1&6&4&1&&&&&2\\
                \hline
                $-3$&&&&&&&3&4&1&&&&\\
                \hline
                $-1$&&&&&&&1&4&2&&&&1\\
                \hline
                $1$&&&&&&&&&3&2&1&&$-2$\\
                \hline
                $3$&&&&&&&&&&1&&&1\\
                \hline
                $5$&&&&&&&&&&&1&1&\\
                \hline
            \end{tabular}
            \caption{Khovanov Polynomial for \texttt{gnoqKDjIMrpEaHblfc}}
        \end{table}
        \begin{table}
            \centering
            \begin{tabular}{| c | c | c | c | c | c | c | c | c | c |}
                \hline
                $q\symbol{92}t$&$-7$&$-6$&$-5$&$-4$&$-3$&$-2$&$-1$&$0$&$\chi$\\
                \hline
                $-19$&2&&&&&&&&$-2$\\
                \hline
                $-17$&&1&&&&&&&1\\
                \hline
                $-15$&&2&3&&&&&&$-1$\\
                \hline
                $-13$&&&1&2&&&&&1\\
                \hline
                $-11$&&&&3&2&&&&1\\
                \hline
                $-9$&&&&&2&2&&&\\
                \hline
                $-7$&&&&&&2&&&2\\
                \hline
                $-5$&&&&&&&2&1&$-1$\\
                \hline
                $-3$&&&&&&&&1&1\\
                \hline
            \end{tabular}
            \caption{Khovanov Polynomial for $9_{49}$}
        \end{table}
        \begin{table}
            \centering
            \begin{tabular}{| c | c |c | c | c | c | c | c | c | c | c | c |}
                \hline
                $q\symbol{92}t$&$-9$&$-8$&$-7$&$-6$&$-5$&$-4$&$-3$&$-2$&$-1$&$0$&$\chi$\\
                \hline
                $-19$&2&&&&&&&&&&$-2$\\
                \hline
                $-17$&&1&&&&&&&&&1\\
                \hline
                $-15$&&2&3&&&&&&&&$-1$\\
                \hline
                $-13$&&&1&3&1&&&&&&1\\
                \hline
                $-11$&&&&3&2&&&&&&1\\
                \hline
                $-9$&&&&&3&4&1&&&&\\
                \hline
                $-7$&&&&&&2&1&1&&&2\\
                \hline
                $-5$&&&&&&&3&2&&&$-1$\\
                \hline
                $-3$&&&&&&&&1&1&1&1\\
                \hline
                $-1$&&&&&&&&&1&1&\\
                \hline
            \end{tabular}
            \caption{Khovanov Polynomial for \texttt{lFKJIOAEnDCpBhmG}}
        \end{table}
        \begin{table}
            \centering
            \begin{tabular}{| c | c | c | c | c | c | c | c | c | c |}
                \hline
                $q\symbol{92}t$&$0$&$1$&$2$&$3$&$4$&$5$&$6$&$7$&$\chi$\\
                \hline
                $5$&1&&&&&&&&1\\
                \hline
                $7$&1&1&&&&&&&\\
                \hline
                $9$&&&1&&&&&&1\\
                \hline
                $11$&&&1&1&1&&&&1\\
                \hline
                $13$&&&&1&2&&&&1\\
                \hline
                $15$&&&&&1&1&&&\\
                \hline
                $17$&&&&&&2&1&&$-1$\\
                \hline
                $19$&&&&&&&&&\\
                \hline
                $21$&&&&&&&&1&$-1$\\
                \hline
            \end{tabular}
            \caption{Khovanov Polynomial for $10_{128}$}
        \end{table}
        \begin{table}
            \centering
            \begin{tabular}{| c | c | c | c | c | c | c | c | c | c | c | c | c |}
                \hline
                $q\symbol{92}t$&$-1$&$0$&$1$&$2$&$3$&$4$&$5$&$6$&$7$&$8$&$9$&$\chi$\\
                \hline
                $1$&1&1&&&&&&&&&&\\
                \hline
                $3$&&1&1&&&&&&&&&\\
                \hline
                $5$&&1&1&1&&&&&&&&1\\
                \hline
                $7$&&&&2&2&&&&&&&\\
                \hline
                $9$&&&&1&1&1&&&&&&1\\
                \hline
                $11$&&&&&1&3&1&&&&&1\\
                \hline
                $13$&&&&&&&1&2&&&&1\\
                \hline
                $15$&&&&&&&1&1&&&&\\
                \hline
                $17$&&&&&&&&&2&1&&$-1$\\
                \hline
                $19$&&&&&&&&&&&&\\
                \hline
                $21$&&&&&&&&&&&1&$-1$\\
                \hline
            \end{tabular}
            \caption{Khovanov Polynomial for \texttt{eHPNqGJlBFoiaDCkM}}
        \end{table}
        \begin{table}
            \centering
            \begin{tabular}{| c | c | c | c | c | c | c | c | c | c | c | c |}
                \hline
                $q\symbol{92}t$&$0$&$1$&$2$&$3$&$4$&$5$&$6$&$7$&$8$&$9$&$\chi$\\
                \hline
                $3$&1&1&&&&&&&&&\\
                \hline
                $5$&1&&&&&&&&&&1\\
                \hline
                $7$&&&2&2&&&&&&&\\
                \hline
                $9$&&&&&1&&&&&&1\\
                \hline
                $11$&&&&1&3&1&&&&&1\\
                \hline
                $13$&&&&&&1&2&&&&1\\
                \hline
                $15$&&&&&&1&1&&&&\\
                \hline
                $17$&&&&&&&&2&1&&$-1$\\
                \hline
                $19$&&&&&&&&&&&\\
                \hline
                $21$&&&&&&&&&&1&$-1$\\
                \hline
            \end{tabular}
            \caption{Khovanov Polynomial for \texttt{edjkaGIlFbch}}
        \end{table}
        \begin{table}
            \centering
            \begin{tabular}{| c | c | c | c | c | c | c | c | c | c |}
                \hline
                $q\symbol{92}t$&$-3$&$-2$&$-1$&$0$&$1$&$2$&$3$&$4$&$\chi$\\
                \hline
                $-9$&1&&&&&&&&$-1$\\
                \hline
                $-7$&&1&&&&&&&1\\
                \hline
                $-5$&&1&1&&&&&&\\
                \hline
                $-3$&&&1&2&&&&&1\\
                \hline
                $-1$&&&&2&1&&&&1\\
                \hline
                $1$&&&&1&2&1&&&\\
                \hline
                $3$&&&&&&1&1&&\\
                \hline
                $5$&&&&&&&1&&$-1$\\
                \hline
                $7$&&&&&&&&1&1\\
                \hline
            \end{tabular}
            \caption{Khovanov Polynomial for $10_{136}$}
        \end{table}
        \begin{table}
            \centering
            \begin{tabular}{| c | c | c | c | c | c | c | c | c | c |}
                \hline
                $q\symbol{92}t$&$-5$&$-4$&$-3$&$-2$&$-1$&$0$&$1$&$2$&$\chi$\\
                \hline
                $-9$&1&&&&&&&&$-1$\\
                \hline
                $-7$&&1&&&&&&&1\\
                \hline
                $-5$&&1&1&&&&&&\\
                \hline
                $-3$&&&1&2&&&&&1\\
                \hline
                $-1$&&&&1&1&1&&&1\\
                \hline
                $1$&&&&&2&2&&&\\
                \hline
                $3$&&&&&&1&1&&\\
                \hline
                $5$&&&&&&&1&&$-1$\\
                \hline
                $7$&&&&&&&&1&1\\
                \hline
            \end{tabular}
            \caption{Khovanov Polynomial for \texttt{igDKHJaEbFC}}
        \end{table}
        \begin{table}
            \centering
            \begin{tabular}{| c | c | c | c | c | c | c | c | c | c | c | c |}
                \hline
                $q\symbol{92}t$&$-9$&$-8$&$-7$&$-6$&$-5$&$-4$&$-3$&$-2$&$-1$&$0$&$\chi$\\
                \hline
                $-21$&1&&&&&&&&&&$-1$\\
                \hline
                $-19$&&&&&&&&&&&\\
                \hline
                $-17$&&1&1&&&&&&&&\\
                \hline
                $-15$&&&&1&1&&&&&&\\
                \hline
                $-13$&&&&1&&&&&&&1\\
                \hline
                $-11$&&&&&1&2&1&&&&\\
                \hline
                $-9$&&&&&&&&&&&\\
                \hline
                $-7$&&&&&&&1&1&&&\\
                \hline
                $-5$&&&&&&&&&&1&1\\
                \hline
                $-3$&&&&&&&&&&1&1\\
                \hline
            \end{tabular}
            \caption{Khovanov Polynomial for $10_{145}$}
        \end{table}
        \begin{table}
            \centering
            \begin{tabular}{| c | c | c | c | c | c | c | c | c | c | c | c | c |}
                \hline
                $q\symbol{92}t$&$-9$&$-8$&$-7$&$-6$&$-5$&$-4$&$-3$&$-2$&$-1$&$0$&$1$&$\chi$\\
                \hline
                $-21$&1&&&&&&&&&&&$-1$\\
                \hline
                $-19$&&&&&&&&&&&&\\
                \hline
                $-17$&&1&1&&&&&&&&&\\
                \hline
                $-15$&&&&1&1&&&&&&&\\
                \hline
                $-13$&&&&1&&&&&&&&1\\
                \hline
                $-11$&&&&&1&2&1&&&&&\\
                \hline
                $-9$&&&&&&&1&1&&&&\\
                \hline
                $-7$&&&&&&&1&1&&&&\\
                \hline
                $-5$&&&&&&&&1&1&1&&1\\
                \hline
                $-3$&&&&&&&&&&1&&1\\
                \hline
                $-1$&&&&&&&&&&1&1&\\
                \hline
            \end{tabular}
            \caption{Khovanov Polynomial for \texttt{eoHKqGJnCFmPDibaL}}
        \end{table}
        \begin{table}
            \centering
            \begin{tabular}{| c | c | c | c | c | c | c | c | c | c | c | c | c | c |}
                \hline
                $q\symbol{92}t$&$-11$&$-10$&$-9$&$-8$&$-7$&$-6$&$-5$&$-4$&$-3$&$-2$&$-1$&$0$&$\chi$\\
                \hline
                $-21$&1&&&&&&&&&&&&$-1$\\
                \hline
                $-19$&&&&&&&&&&&&&\\
                \hline
                $-17$&&1&1&&&&&&&&&&\\
                \hline
                $-15$&&&&1&1&&&&&&&&\\
                \hline
                $-13$&&&&1&&&&&&&&&1\\
                \hline
                $-11$&&&&&1&2&1&&&&&&\\
                \hline
                $-9$&&&&&&&&1&1&&&&\\
                \hline
                $-7$&&&&&&&1&1&&&&&\\
                \hline
                $-5$&&&&&&&&&1&2&&&1\\
                \hline
                $-3$&&&&&&&&&&&&1&1\\
                \hline
                $-1$&&&&&&&&&&&1&1&\\
                \hline
            \end{tabular}
            \caption{Khovanov Polynomial for \texttt{kNJIpHLFECoMGABd}}
        \end{table}
        \begin{table}
            \centering
            \begin{tabular}{| c | c | c | c | c | c | c | c | c | c | c | c |}
                \hline
                $q\symbol{92}t$&$-9$&$-8$&$-7$&$-6$&$-5$&$-4$&$-3$&$-2$&$-1$&$0$&$\chi$\\
                \hline
                $-23$&1&&&&&&&&&&$-1$\\
                \hline
                $-21$&&&&&&&&&&&\\
                \hline
                $-19$&&1&1&&&&&&&&\\
                \hline
                $-17$&&&&1&1&&&&&&\\
                \hline
                $-15$&&&&1&1&&&&&&\\
                \hline
                $-13$&&&&&1&2&1&&&&\\
                \hline
                $-11$&&&&&&1&&&&&1\\
                \hline
                $-9$&&&&&&&1&1&&&\\
                \hline
                $-7$&&&&&&&&&&1&1\\
                \hline
                $-5$&&&&&&&&&&1&1\\
                \hline
            \end{tabular}
            \caption{Khovanov Polynomial for $10_{161}$}
        \end{table}
        \begin{table}
            \centering
            \begin{tabular}{| c | c | c | c | c | c | c | c | c | c | c | c | c |}
                \hline
                $q\symbol{92}t$&$-9$&$-8$&$-7$&$-6$&$-5$&$-4$&$-3$&$-2$&$-1$&$0$&$1$&$\chi$\\
                \hline
                $-23$&1&&&&&&&&&&&$-1$\\
                \hline
                $-21$&&&&&&&&&&&&\\
                \hline
                $-19$&&1&1&&&&&&&&&\\
                \hline
                $-17$&&&&1&1&&&&&&&\\
                \hline
                $-15$&&&&1&1&&&&&&&\\
                \hline
                $-13$&&&&&1&2&1&&&&&\\
                \hline
                $-11$&&&&&&1&1&1&&&&1\\
                \hline
                $-9$&&&&&&&1&1&&&&\\
                \hline
                $-7$&&&&&&&&1&1&1&&1\\
                \hline
                $-5$&&&&&&&&&&1&&1\\
                \hline
                $-3$&&&&&&&&&&1&1&\\
                \hline
            \end{tabular}
            \caption{Khovanov Polynomial for \texttt{hOqrljsnMeipFAgkbcd}}
        \end{table}
        \begin{table}
            \centering
            \begin{tabular}{| c | c | c | c | c | c | c | c | c | c |}
                \hline
                $q\symbol{92}t$&$-7$&$-6$&$-5$&$-4$&$-3$&$-2$&$-1$&$0$&$\chi$\\
                \hline
                $-19$&1&&&&&&&&$-1$\\
                \hline
                $-17$&&1&&&&&&&1\\
                \hline
                $-15$&&1&2&&&&&&$-1$\\
                \hline
                $-13$&&&1&2&&&&&1\\
                \hline
                $-11$&&&&2&2&&&&\\
                \hline
                $-9$&&&&&2&1&&&$-1$\\
                \hline
                $-7$&&&&&&2&1&&1\\
                \hline
                $-5$&&&&&&&1&1&\\
                \hline
                $-3$&&&&&&&&2&2\\
                \hline
            \end{tabular}
            \caption{Khovanov Polynomial for \texttt{K11n118}}
        \end{table}
        \begin{table}
            \centering
            \begin{tabular}{| c | c | c | c | c | c | c | c | c | c | c | c | c | c | c | c |}
                \hline
                $q\symbol{92}t$&$-7$&$-6$&$-5$&$-4$&$-3$&$-2$&$-1$&$0$&$1$&$2$&$3$&$4$&$5$&$\chi$\\
                \hline
                $-19$&1&&&&&&&&&&&&&$-1$\\
                \hline
                $-17$&&1&&&&&&&&&&&&1\\
                \hline
                $-15$&&1&2&&&&&&&&&&&$-1$\\
                \hline
                $-13$&&&1&3&1&&&&&&&&&1\\
                \hline
                $-11$&&&&2&2&&&&&&&&&\\
                \hline
                $-9$&&&&&4&5&2&&&&&&&$-1$\\
                \hline
                $-7$&&&&&&2&2&1&&&&&&1\\
                \hline
                $-5$&&&&&&1&4&4&1&&&&&\\
                \hline
                $-3$&&&&&&&&3&3&2&&&&2\\
                \hline
                $-1$&&&&&&&&1&2&1&&&&\\
                \hline
                $1$&&&&&&&&&&2&3&1&&\\
                \hline
                $3$&&&&&&&&&&&&&&\\
                \hline
                $5$&&&&&&&&&&&&1&1&\\
                \hline
            \end{tabular}
            \caption{Khovanov Polynomial for \texttt{bfINqJaKLTGHPRsCeMDo}}
        \end{table}
    \newpage
    \bibliographystyle{plain}
    \bibliography{bib.bib}
    \newpage
    The source code used to generate this document, including figures, and all
    of the analysis is free software and released under version 3 of the
    GNU General Public License. See
    \cite{MaguireMathematicsAndPhysics} and \cite{MaguireLibtmpl}.
    \par\hfill\par
    Vladimir Chernov
    \par
    6188 Kemeny Hall
    \par
    Mathematics Department, Dartmouth College
    \par
    Hanover, NH USA 03755
    \par
    Vladimir.Chernov@dartmouth.edu
    \par\hfill\par
    Ryan Maguire
    \par
    2-106
    \par
    Mathematics Department, MIT
    \par
    Cambridge, MA USA 02139
    \par
    rmaguire@mit.edu
\end{document}